\documentclass[11pt]{article}
\usepackage{epsfig,pstricks,exscale,here}
\usepackage{fullpage}
\usepackage{amsthm,amsmath,enumerate}
\usepackage{latexsym}
\usepackage{amssymb,bbm}
\usepackage{wrapfig}

\newcommand{\N}{\mathbb{N}} 
\newcommand{\R}{\mathbb{R}} 
\newcommand{\C}{\mathbb{C}} 
\newcommand{\SH}{\mathbb{S}}
\newcommand{\Proj}{\mathbb{CP}}

\newtheorem{theorem}{Theorem}
\newtheorem{lemma}[theorem]{Lemma}
\newtheorem{propo}[theorem]{Proposition}
\newtheorem{coro}[theorem]{Corollary}
\newtheorem{algorithm}[theorem]{Algorithm}
\newtheorem{defini}[theorem]{Definition}
\newtheorem{conj}[theorem]{Conjecture}

\newtheorem{inv}[theorem]{Invariant}

\def\rc{\advance\leftskip by 0pt plus 40em\rightskip=\leftskip
  \parfillskip=0pt \spaceskip=.3333em \xspaceskip=.5em
  \pretolerance=9999 \tolerance=9999 \hyphenpenalty=9999
  \exhyphenpenalty=9999 }

\title{The Szemer\'{e}di-Trotter Theorem in the  Complex Plane}
\author{Csaba D. T\'oth\thanks{Department of Mathematics, California State University Northridge, Los Angeles, CA; and Department of Computer Science, Tufts University, Medford, MA, USA. Email: \texttt{cdtoth@acm.org}
Research on this paper was conducted at the E\"otv\"os Lor\'and University, Budapest.}}
\date{}
\begin{document}

\maketitle

\vspace{-\baselineskip}

\begin{abstract}
It is shown that $n$ points and $e$ lines in the complex Euclidean
plane $\C^2$ determine $O(n^{2/3}e^{2/3}+n+e)$ point-line incidences.
This bound is the best possible, and it generalizes the celebrated
theorem by Szemer\'edi and Trotter about point-line incidences in the
real Euclidean plane $\R^2$.
\end{abstract}

\section{Introduction}

It was shown by Szemer\'edi and Trotter~\cite{SzT:83}, settling a conjecture by Erd\H{o}s,
that $n$ points and $e$ lines in the Euclidean plane $\R^2$ determine at most
$O(n^{2/3}e^{2/3}+n+e)$ point-line incidences. This bound is the best possible~\cite{Ed87,E82}:
there are $\Theta(n^{2/3}e^{2/3}+n+e)$ point-line incidences between
the points of an $\lfloor\sqrt{n}\rfloor \times \lfloor\sqrt{n}\rfloor$ section of the
integer lattice and $e$ appropriate lines in $\R^2$. Originally,
Szemer\'edi and Trotter~\cite{SzT:83} proved an upper bound of $10^{60}n^{2/3}e^{2/3}+3n+3e$
for the number of point-line incidences in $\R^2$. The constant coefficients have been 
improved substantially, and the current best upper bound~\cite{PRTT06} is $2.5n^{2/3}e^{2/3}+n+m$.

The Szemer\'edi-Trotter bound is fundamental in combinatorial geometry~\cite{PS:04}. Its importance is illustrated by the fact that, since the original publication of their result, two completely different proof techniques have been developed for it, each of which has many applications on its own right. One is the theory of \emph{$\varepsilon$-cuttings} and \emph{$\varepsilon$-nets} based on a divide-and-conquer strategy~\cite{CEGSW:90}, the other is the \emph{crossing number theory} for graphs drawn in the plane~\cite{PS:96,Sz:96}. The Szemer\'edi-Trotter bound has several generalizations to point-curve incidences in $\R^2$, where the curves are pseudo-lines~\cite{CEGSW:90}, bounded degree algebraic curves~\cite{PS:92}, or Jordan curves with certain intersection constraints~\cite{PS:96}.

Extending some of the applications of the Szemer\'edi-Trotter bound requires
a similar bound for points and lines in the \emph{complex} Euclidean plane $\mathbb{C}^2$.
However, all existing proofs of the Szemer\'edi-Trotter theorem heavily rely
on the topology of the \emph{real} Euclidean plane $\mathbb{R}^2$ and no natural complex or
multidimensional counterparts have been found so far. A line in the complex
plane is the set of points $(x,y)\in \C^2$ satisfying a linear equation
$y=ax+b$ or $x=b$ for some $a,b\in \C$. Our main result is the following.

\begin{theorem}\label{main}
There is a constant $C$ such that $n$ points and $e$ lines in the
complex Euclidean plane $\C^2$ determine at most
$C n^{2/3}e^{2/3} + 3n + 3e$ point-line incidences.
\end{theorem}

The upper bound in Theorem~\ref{main} is the best possible apart from constant factors.
Consider a configuration of $n$ points and $e$ lines in $\R^2$ that
determine $\Theta(n^{2/3}e^{2/3} + n + e)$ point-line incidences~\cite{Ed87}.
Every point $(a,b)\in \R^2$ can be embedded as a point $(a,b)\in \C^2$, 
and every line $y=cx+d$ or $x=d$, with $c,d\in \R$, is contained in the complex
line $y=cx+d$ or $y=d$, with $c,d\in \C^2$, with the same point-line incidences.

The proof of Theorem~\ref{main} is presented in Sections~\ref{pre}--\ref{CombL}.
It is essentially the adaptation of the original proof by Szemer\'edi and Trotter to the complex plane.
We prove that the constant $C$ in Theorem~\ref{main} may be taken to be $C=10^{60}$.
No effort has been made to optimize the value of $C$. In order to keep the presentation
as simple as possible, the constants are often estimated very generously. We note here,
though, that the proof technique inevitably leads to a large constant $C$, similarly to
the original proof by Szemer\'edi and Trotter~\cite{SzT:83}.

\paragraph{Corollaries.}
We present a few immediate consequences of Theorem~\ref{main}.
The analogues of these results in the real Euclidean plane can be derived
from the Szemer\'edi-Trotter theorem by purely combinatorial arguments,
so they immediately generalize to the complex Euclidean plane.
An equivalent formulation of the Szemer\'edi-Trotter theorem is an
upper bound on the number of lines containing at least $t$, $2\leq
t\leq n$, points in the plane~\cite{SzT:83}.

\begin{theorem}\label{eqq}
For $n$ points in $\C^2$, and an integer $t$, $2\leq t\leq n$,
the number of complex lines incident to at least $t$ points is
$$O\left( \frac{n^2}{t^3} + \frac{n}{t} \right).$$
\end{theorem}

A result by Beck~\cite{B:83} follows from the Szemer\'edi-Trotter theorem
by purely combinatorial arguments, and hence it generalizes to the complex plane.

\begin{coro}
There is a constant $c_1>0$ such that, for every set of $n$ points in $\C^2$,
at least one of the following two statements holds.
\begin{itemize}\itemsep -2pt
\item There are at least $c_1n^2$ complex lines, each of which is incident to at least
      two points.
\item There is a complex line incident to at least $n/100$ points.
\end{itemize}
\end{coro}

For a set $A\subset \C$, the set of pairwise sums and products formed by the elements of
$A$ is $A+A=\{a+b: a,b\in A\}$ and $A\cdot A=\{a\cdot b:a,b\in A\}$, respectively.
Elekes~\cite{El1} proved that for every finite set $A\subset \R$, $\max \{ |A+A|,
|A\cdot A|\}=\Omega (|A|^{5/4})$. This bound was later improved to
$\Omega(|A|^{14/11})$ by Solymosi~\cite{S}. The same combinatorial
argument over complex numbers yields the following.

\begin{coro}
There is a constant $c_2>0$ such that for every
finite set $A\subset \C$, we have
$$c_2\cdot |A|^{14/11} \leq \max \{ |A+A|, |A\cdot A| \}.$$
\end{coro}

Our last corollary generalizes a theorem by Elekes~\cite{El2} from
\emph{homothetic} subsets of $\R^2$ to \emph{similar} subsets of $\R^2$.
For two finite point sets $A,A'\subset \R^2$, we denote by $A\sim A'$ if
they are \emph{similar} to each other (that is, an isometry followed by
a dilation maps $A$ to $A'$). For two finite point sets, $A,B \subset \R^2$,
let
$$S(A,B) = | \{A'\subset \R^2 : A'\subset B, A' \sim A \} |$$
be the number of similar copies of $A$ in $B$. The maximal number of
similar copies of a set of $t$ points in a set of $n$ points
in $\R^2$ is denoted by $s(t,n)=\max\{ S(A,B) : |A|=t, |B|=n \}.$

\begin{coro}
There is a constant $c_3$, such that for every $t,n\in\N$, we have
$$s(t,n)\leq \frac{c_3 n^2}{t}.$$
\end{coro}

\paragraph{Possible generalizations.}
A natural generalization of the Szemer\'edi-Trotter theorem (and our
Theorem~\ref{main}) would be an upper bound on the number of
incidences between points and $d$-flats in $\R^{2d}$.
A \emph{$d$-flat} in Euclidean space is a $d$-dimensional affine subspace.

\begin{conj}\label{cony}
For every integer $d \geq 1$, there is a constant $c_d$ with the following property.
Given $n$ points and $e$ $d$-flats in $\R^{2d}$ such that the intersection of every two $d$-flats
is either empty or a single point, then they determine at most $c_d(n^{2/3}e^{2/3} + n + e)$
incidences.
\end{conj}

For $d=1$, this is equivalent to the Szemer\'edi-Trotter theorem.
Our Theorem~\ref{main} is a special case for $d=2$ where all 2-flats
correspond to complex lines in $\C^2$.

\section{Preliminaries~\label{pre}}

The proof of Theorem~\ref{main} was conceived in an attempt to prove
Conjecture~\ref{cony}, by generalizing the original proof of the Szemer\'edi-Trotter theorem.
Most of the steps of the proof are either purely combinatorial or use the real Euclidean space,
where lines in $\C^2$ are embedded as 2-flats in $\R^4$. The use of special properties
of lines in $\C^2$ (as opposed to 2-flats in $\R^4$) is kept to a minimum.

In Subsection~\ref{ssec:out}, we present a brief outline of the original proof of Szemer\'edi and Trotter~\cite{SzT:83}, and point out the similarities and key differences from our proof. Subsection~\ref{ssec:gr} summarizes the few basic properties of Grassmann manifolds. We exploit
properties of $\C^2$ only in Subsection~\ref{ssec:cmplx}, the proof of our \emph{Separation Lemma}. Subsections~\ref{ssec:comb} and~\ref{ssec:dis} use purely combinatorial arguments, and
Sections~\ref{CovL} and~\ref{CombL} rely exclusively on the geometry of $\R^d$.

\subsection{Outline\label{ssec:out}}

Our proof follows the same strategy as that of Szemer\'edi and Trotter~\cite{SzT:83}.
We briefly outline these steps below.

\begin{itemize}
\item[{\rm (i)}]
The proof by Szemer\'edi and Trotter proceeds by contradiction, and considers
a minimum counterexample, that is, a system $(P,E)$ of
$n$ points and $e$ lines with more than $Cn^{2/3}e^{2/3}+3n+3e$
point-line incidences, where $n+e$ is minimal.
\item[{\rm (ii)}]
They show that $(P,E)$ contains a large subsystem, $(P_3,L_1\cup L_2)$,
that has a ``regular'' structure (Separation Lemma). In particular, $P_3\subset P$
contains $\Omega(n)$ points; each point in $P_3$ is incident to $\Omega(e^{2/3}/n^{1/3})$
lines in each of $L_1$ and $L_2$; and the directions of the lines in $L_1$ and $L_2$ are
close to two orthogonal directions after an appropriate linear transformation.
\item[{\rm (iii)}]
 The Covering Lemma by Szemer\'edi and Trotter~\cite{SzT:83b} shows that for a point set $P_3$ in the plane and a parameter $1\leq k\leq |P_3|$, there are interior-disjoint axis-parallel squares such that each square contains $\Theta(k)$ points of $P_3$ and they jointly contain $\Omega(n)$ points in $P_3$.
\item[{\rm (iv)}]
The combination of the Separation Lemma and the Covering Lemma leads to a lower bound for
the total number of intersection points (crossings) between lines in $L_1$ and $L_2$.
One can find at least $e^2$ crossings between lines in $L_1\subset E$ and $L_2\subset E$,
contradicting the fact that there cannot be more than ${e\choose 2}$ crossings between lines in $E$.
\end{itemize}

We prove a generalization of the Separation Lemma for lines in $\C^2$ in Section~\ref{SepL}.
We prove an elaborate version of the Covering Lemma in Section~\ref{CovL};
and give a lower bound for the number of intersection points of complex lines in Section~\ref{CombL}.

In the complex plane, however, there are several key differences compared to the
real plane. Szemer\'edi and Trotter~\cite{SzT:83} cover a constant fraction of
the points by interior-disjoint axis-aligned squares in $\R^2$. They make use of
the simple but crucial fact that when a square in $\R^2$ is decomposed into four right
triangles along the two diagonals then every two points $p$ and $q$ lying in one
right triangle has the following property: if $p$ (resp., $q$) is incident to two
lines $\ell_1^p$ and $\ell_2^p$ (resp., $\ell_1^q$ and $\ell_2^q$) that are almost
parallel to the two diagonals of the square, then at least one of the intersection
points, $\ell_1^p\cap \ell_2^q$ or $\ell_2^p\cap \ell_1^q$, is in the interior of the square.
This property does not generalize to 2-flats in a decomposition of a 4-dimensional hypercube.
This explains why we need a significantly more involved covering lemma in $\R^4$.

Similarly, difficulties arise if we want to find appropriate regular
structures, like those in our Separation Lemma. Szemer\'edi and
Trotter used the space of directions of lines of the Euclidean plane
and found a linear transformation that produces two almost orthogonal
families of lines. The space of directions of lines in $\C^2$ is
much more difficult to handle.

\subsection{Grassmann manifolds\label{ssec:gr}}

A line $\ell$ in the complex plane $\C^2$ is defined by a linear equation $y=ax+b$ or $x=b$ for some $a,b\in \C$. The direction of $\ell$ can be represented by the parallel line incident to the origin, $y=ax$ or $x=0$, respectively, or by its \emph{slope}, which is $a$ or $\infty$, respectively.
The space of 1-dimensional subspaces in $\C^2$ is the Grassmann manifold $H(1,1)$.
It can be represented by the complex projective line $\Proj^1$ or the Riemann sphere $\mathbb{C}\cup \{\infty\}$~\cite{H}.

\begin{figure}
  \centering \vspace{-\baselineskip}
  \includegraphics[width=.4\textwidth]{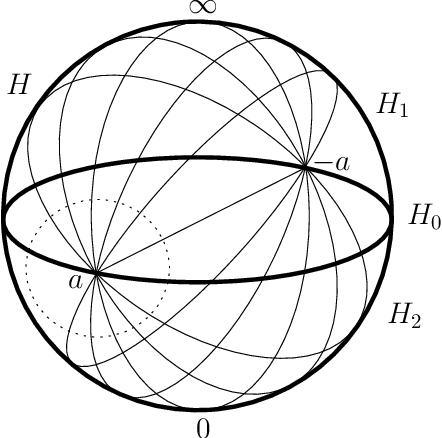}
  \caption{A Riemann sphere $H$, a unit slope $a\in H_0$, its antipodal slope
  $-1/\overline{a}=-a\in H_0$, and some orbits of $\pi_\lambda^a$ between $a$ and $-a$.
  \label{riemann}}
\end{figure}

The standard correspondence between the Riemann sphere $\mathbb{C}\cup \{\infty\}$ and the unit sphere $\SH^2$ is defined as follows. Identify every slope $a\in \C$ with the point $({\rm Re}(a),{\rm Im}(a),0)\in \R^3$ in the plane $z=0$ of $\R^3$. A stereographic projection maps every point from the plane $z=0$ to the unit sphere $H=\{(x,y,z)\in \mathbb{R}^3:x^2+y^2+z^2=1\}$; and the slope $\infty$ is mapped to the ``North Pole'' $(0,0,1)\in H$. Note that every unit slope $a$, $|a|=1$, is mapped to the \emph{equator} $H_0=\{(x,y,z)\in \mathbb{R}^3: x^2+y^2=1, z=0\}$ of $H$~\cite{H}. Denote the two closed hemispheres of $H$ above and below $H_0$ by $H_1$ and $H_2$, respectively (see Fig.~\ref{riemann}).

For a line $\ell$ in $\C^2$, let $\hat{\ell} \in H$ denote its direction.
Similarly, for a set $L$ of complex lines, let $\widehat{L}\subseteq H$ denote
the multiset of directions of the lines in $L$.

\paragraph{Metrics.}
$H(1,1)$ has an essentially unique metric, invariant under unitary transformations.
The distance ${\rm dist}(\hat{\ell}_1,\hat{\ell}_2)$ between the directions
of two lines $\ell_1$ and $\ell_2$ in $\C^2$ can be defined in
terms of their {\em principal angle}
$\rm{arccos} ( \max\{ |\mathbf{u} \overline{\mathbf{v}} | : {\mathbf u}
\in \ell_1, \mathbf{v} \in \ell_2, |\mathbf{u}|=|{\mathbf v}|=1\})$.
This is equivalent to the chordal distance in the Riemann sphere
representation~\cite{GL,W}. In this paper, we always use the chordal
metric in $\SH^2$, measured in degrees. For example, if two directions
$a_1,a_2\in \C\cup \{\infty\}$ are \emph{orthogonal} (i.e.,
$a_1\overline{a}_2=-1$, or $a_1=0$ and $a_2=\infty$, or $a_1=\infty$ and $a_2=0$),
then they correspond to antipodal points in $\SH^2$, hence ${\rm dist}(a_1,a_2)=180^{\circ}$.

\paragraph{Embedding into ${\rm Gr}(2,2)$.} The map $\tau: \C^2\rightarrow \R^4$, $(z_1,z_2)\rightarrow
({\rm Re}(z_1), {\rm Im}(z_1), {\rm Re}(z_2,), {\rm Im}(z_2))$ identifies the points of $\C^2$ and $\R^4$,
it maps the lines in $\C^2$ into 2-flats in $\R^4$. Since parallel lines are mapped to
parallel 2-flats, it induces an embedding $\hat{\tau}$ of $H(1,1)$
into the Grassmann manifold ${\rm Gr}(2,2)$ of 2-flats in $\R^4$.
${\rm Gr}(2,2)$ has several different metrics, invariant under orthogonal
transformations. All metrics can be defined in terms of the two principal
angles between two 2-flats~\cite{H}. We consider the distance between the
directions of two 2-flats in $\R^4$ to be the sum of their principal angles.
In particular, the distance of two orthogonal 2-subspaces is $2\cdot 90^{\circ}=180^{\circ}$.
With this choice, the embedding $\hat{\tau}:H(1,1)\rightarrow {\rm Gr}(2,2)$ preserves
the metric of $H(1,1)$. Specifically, $\hat{\tau}$ maps a $3^\circ$-neighborhood
in $H(1,1)$ into a $3^\circ$-neighborhood in ${\rm Gr}(2,2)$.

\paragraph{Nondegenerate complex linear transformations.}
The group $GL(2,\C)$ of nondegenerate complex linear transformations
acts on $\C^2$ and preserves point-line incidences. The group
$PGL(2,\C)$ acts on $H(1,1)=\Proj^1$, and corresponds to
M\"obius transformations. Each linear transformation in $GL(2,\C)$
induces a M\"obius transformation in $PGL(2,\C)$.
Let $\Psi \subset GL(2,\C)$ denote the set of nondegenerate linear transformations
on $\C^2$ that induce automorphisms on each of $H_0$, $H_1\setminus H_0$, and
$H_2\setminus H_0$.

For every $\lambda \in \R$, $-1<\lambda <1$, we define the linear
transformation $\pi_\lambda^1\in GL(2,\C)$ to be
$$\pi_{\lambda}^1: \C^2\longrightarrow \C^2,\hspace{5mm}
\left[ \begin{array}{c}z_1\\z_2\end{array}\right]\longrightarrow
\frac{1}{\sqrt{1-\lambda^2}}\ \left[ \begin{array}{cc} 1 &
\lambda\\ \lambda & 1 \end{array}\right]
\left[\begin{array}{c}z_1\\z_2\end{array}\right].$$

For every vector $(z_1,z_2)\in \C^2$, the transformation
$\pi_\lambda^1$ dilates the component parallel to $(1,1)\in \C^2$ by
a factor of $\sqrt{\frac{1+\lambda}{1-\lambda}} \in \R$ and the perpendicular
component by a factor of $\sqrt{\frac{1-\lambda}{1+\lambda}}\in \R$. Note that
$\pi_\lambda^1\in \Psi$, since a vector $(z_1,z_2)\in \C^2$ has unit slope
(that is, $|z_1|=|z_2|$) if and only if $\pi_\lambda^1((z_1,z_2))$
has unit slope. The orbit of an element of $H\setminus \{1,-1\}$ under the transformations
$\pi_{\lambda}^1$, $\lambda\in (-1,1)$, is a main halfcircle between $1\in H_0$
and its antipodal $-1 \in H_0$.

For every $a\in \C$, $|a|=1$, we define the linear
transformation $\pi_{\lambda}^a := \varrho^a \pi_{\lambda}^1 (\varrho^a)^{-1}$, where
$$\varrho^a: \C^2 \longrightarrow \C^2,\hspace{5mm}
\left[ \begin{array}{c}z_1\\z_2\end{array}\right]\longrightarrow
\left[ \begin{array}{cc} 1 & 0\\ 0&
1/a\end{array} \right]
\left[
\begin{array}{c}z_1\\z_2\end{array}\right].$$
$\varrho^a$ is a unitary transformation (it induces an isometry on the
Riemann sphere $H$) and $\varrho^a\in \Psi$. The orbits of the elements of
$H\setminus \{a,-a\}$ under $\pi_{\lambda}^a$, $\lambda\in (-1,1)$,
are the main halfcircles between $a\in H_0$ and its antipodal $1/\overline{a}=-a\in H_0$
(see Fig.~\ref{riemann}). As $\lambda$ continuously increases in the interval $(-1,1)$,
the image of every point in $H$ (except for $a$ and $-a$) moves continuously
toward $a\in H_0$ along a main halfcircle of $H$.

We use one more property of the space of complex directions in the
proof of our Separation Lemma. Two sufficiently small disjoint disks
(i.e., spherical caps) in the Riemann sphere $H$ can be mapped into
two small neighborhoods around two orthogonal directions by a nondegenerate
linear transformation of $\mathbb{C}^2$ if the disks are at least a constant 
distance apart. In contrast, ${\rm Gr}(2,2)$ does not have this property: 
if two disjoint disks in ${\rm Gr}(2,2)$ contain 2-subspaces whose 
intersection is a line, then no linear transformation of $\mathbb{R}^4$ 
can increase their distance above $90^{\circ}$. This is why our proof 
technique cannot establish Conjecture~\ref{cony} for $d=2$ (that is, 
complex lines cannot be replaced by arbitrary 2-flats in $\R^4)$.

\section{Separation Lemma \label{SepL}}

Our Separation Lemma (Lemma~\ref{first}) is a straightforward
generalization of a corresponding result by Szemer\'edi and Trotter in $\R^2$.
It claims that a minimal counterexample to Theorem~\ref{main} contains
a large and fairly regular structure of points and lines.
Let $(P,E)$ be a \emph{system} of a finite set of points $P$ and a finite
set of lines $E$ in $\C^2$. Let $n=|P|$ and $e=|E|$ denote the number of
points and lines, respectively, and let $I=I_{P,E}$ denote the number
of point-line incidences between $P$ and $E$. A system $(P,E)$ is
\emph{critical} if $e+n$ is minimal among all systems where
the number of incidences satisfies $I>\max(C n^{2/3}e^{2/3}, 3n ,3e)$
with constant $C=10^{60}$.

\begin{lemma} {\bf (Separation Lemma)}\label{first}
If $(P,E)$ is a critical system, then there is a point set
$P_3 \subseteq P$ and two line sets $L_1,L_2\subset E$ such
that for the constant $M=10^{10}$, we have
\begin{itemize}\itemsep -2pt
\item[{\rm (a)}] $|P_3| \geq n/M^8$;
\item[{\rm (b)}] every point $p\in P_3$ is incident to at least $I/(nM^2)$ lines from each of $L_1$ and $L_2$;
\item[{\rm (c)}] there are two orthogonal directions $\hat{\ell}_1, \hat{\ell}_2\in H(1,1)$ such that an appropriate nondegenerate linear transformation of $\C^2$ maps the directions of the lines in $L_1$ and $L_2$ into the $3^\circ$-neighborhood of $\hat{\ell}_1$ and $\hat{\ell}_2$, respectively.
\end{itemize}
\end{lemma}

\subsection{Combinatorial preprocessing\label{ssec:comb}}

First we show that in a critical system $(P,E)$, the parameters
$n$ and $e$ cannot be too far from each other, more precisely,
each of them is much larger than the square root of the other.

\begin{lemma} \label{alap}
If $(P,E)$ is a critical system, then
\begin{equation}\label{1}
e > \frac{C^3}{3^{3/2}} \sqrt{n},
\hspace{1cm} \mbox{\rm and}\hspace{1cm}
n > \frac{C^3}{3^{3/2}}\sqrt{e}.
\end{equation}
\end{lemma}

\proof
By symmetry, it is enough to prove the first inequality. For every point
$p\in P$, denote by  $d_p$ the number of lines in $E$ incident to $p$.
By Jensen's inequality, we have
\begin{eqnarray}
e^2 &>& {e\choose 2 } \geq \sum_{p\in P}{ d_p\choose 2 } =
\sum_{p\in P} \frac{d_p^{2}}{2} - \sum_{p\in P} \frac{d_p}{2}
\geq  \frac{1}{2n} \left(\sum_{p\in P}  d_p\right)^2 - \frac{1}{2}
\sum_{p\in P}  d_p\nonumber\\
&=& \frac{1}{2n} I^{2} - \frac{1}{2} I> \frac{1}{2n}
I^{2}-\frac{1}{6n} I^2 =\frac{I^2}{3n},\nonumber
\end{eqnarray}
where the last step follows from $I> 3n$.  Therefore, by $I> C
n^{2/3}e^{2/3}$, we have $$e^2> \frac{C^2}{3}n^{1/3}e^{4/3}. \hspace{6mm} \qed$$

\begin{coro}\label{en} If $({P},{E})$ is critical, then we have
$$e=e^{1/3}e^{2/3}< \frac{3}{C^2} n^{2/3}e^{2/3}
<\frac{3}{C^3}I \hspace{3mm} {\rm and} \hspace{3mm} n=n^{1/3}n^{2/3}<
\frac{3}{C^2} n^{2/3}e^{2/3}<\frac{3}{C^3}I. \hspace{6mm} \qed$$
\end{coro}

\begin{coro}\label{enne}
If $(P,E)$ is critical, then
$$\max(C n^{2/3}e^{2/3},3n,3e)=C n^{2/3}e^{2/3}.\qed$$
\end{coro}

Our next goal is to show that each point in $P$ is incident to a large
number of lines in $E$.
Let $d_A=I/n$ denote the average number of lines in $E$ incident to a point in $P$,
and let $f_A=I/e$ be the average number of points in $P$ incident to a line in $E$.
For a point $p\in \C^2$ and a set $F$ of lines in $\C^2$,
denote by $F^p\subseteq F$ the subset of lines in $F$ incident to $p$.
We show that every point in $P$ is incident to at least half the average number of lines.

\begin{lemma}\label{fele}
If $(P,E)$ is critical, then every point in $P$ is incident to at least
$d_A/2$ lines of $E$ and every line in $E$ is incident to at least $f_A/2$ points of $P$.
\end{lemma}

\proof By symmetry, it is enough to prove the first claim,
that is, $|E^p| \geq d_A/2$ for every $p\in P$.
Suppose, to the contrary, that there is a point $p\in P$ incident
to fewer than $d_A/2$ lines in $E$.

Since the system $(P\setminus \{p\},E)$ is smaller than the critical system
$(P,E)$, we have $I_{P\setminus \{p \},E} \leq \max( C
(n-1)^{2/3}e^{2/3},3(n-1),3e)$. This, together with
Corollary~\ref{en}, implies an upper bound on the total number of
incidences in the system $(P,E)$.
$$ I < \frac{d_A}{2} + \max(C (n-1)^{2/3} e^{2/3},3(n-1),3e) <
\frac{1}{2n} \cdot I+ \max \left (\left(\frac{n-1}{n}\right)^{2/3},
\frac{10}{C^3}\right)\cdot I,$$
$$ 1 < \frac{1}{2n} + \max \left( \left(\frac{n-1}{n}\right)^{2/3},
\frac{10}{C^3}\right).$$
The last inequality is equivalent to either $4 n^2 -2 n - 1 < 0$ or
$n< 1/(2\cdot (1-10/C^3))$ depending on the value in the
maximum. Neither inequality has any positive integer solution.
\qed

\subsection{Distinguishing two sets of lines\label{ssec:dis}}

Recall that we represent the space of directions of complex lines as the Riemann
sphere $H$, where a main circle $H_0$ corresponds to the directions of unit slope,
and the two closed hemispheres of $H$ bounded by $H_0$ are denoted $H_1$ and $H_2$, respectively.
We may assume, after applying a nondegenerate linear transformation of $\C^2$,
that $|H_1\cap \widehat{E}|\geq e/2$, $|H_2\cap \widehat{E}|\geq e/2$,
and there is at most one family of parallel lines in $E$ whose
direction is in $H_0$.

\begin{defini}
Let $E=E_1\cup E_2$ be a partition of the line set $E$ such that
$E_1\subset \{\ell\in E: \hat{\ell}\in H_1\}$,
$E_2\subset \{\ell \in E: \hat{\ell} \in H_2\}$,
$|E_1|= \lfloor e/2\rfloor$ and $|{E}_2 |=\lceil e/2\rceil$.
\end{defini}

Note that the parallel lines in $E$ whose direction corresponds
to a point in $H_0$ may belong to $E_1$ or $E_2$.

\begin{defini}
Let $P_0=\{ p\in P: |E_1^p|\geq d_A/100$ and
$|E_2^p|\geq d_A/100\}$ be the set of points incident to at least
$d_A/100$ lines from each of $E_1$ and $E_2$. Partition the points of
$P\setminus P_0$ into two subsets: let $P_1= \{ p \in P \setminus
P_0 : |E_1^p| > |E_2^p| \}$ and $P_2= \{ p \in P \setminus P_0 :
|E_1^p| \leq |E_2^p| \}$.
\end{defini}

\begin{lemma}
If $(P,E)$ is critical, then $|P_0| \geq n/10$.
\end{lemma}

\proof For $j=0,1,2$, let $|P_j| = x_j n$ for some $x_j\geq 0$.
Suppose, to the contrary, that $x_0 < \frac{1}{10}$.
Let $I_j$ denote the number of incidences of the system $(P_j,E)$.
Then the total number of incidences is $I=\sum_{j=0}^{2} I_j$.

Note that $P\neq P_1$, otherwise there would be at most $I/100$ incidences
between $P$ and $E_2$, hence some line in $E_2$ would be incident to fewer than
$f_A/2$ points, contradicting Lemma~\ref{fele}. One can show analogously that $P\neq P_2$.
It follows that each of the systems $(P_0, E)$, $(P_1, E_1)$, and $(P_2, E_2)$
is smaller than the critical system $(P,E)$, hence the bound of Theorem~\ref{main}
applies to each of them. Taking into account the incidences of the systems
$(P_1, E_2)$ and $(P_2, E_1)$, as well, we obtain:
\begin{eqnarray*}
I_{0} &<& C(x_{0}n)^{2/3}e^{2/3}  + 3x_{0}n +3e,\\
I_{1} &<& C(x_{1}n)^{2/3} \lfloor e/2\rfloor^{2/3}+3x_{1}n+3\lfloor e/2\rfloor +(x_{1}n)(d_{A}/100),\\
I_{2} &<& C(x_{2}n)^{2/3} \lceil e/2\rceil^{2/3}+3x_{2}n+3\lceil e/2\rceil +(x_{2}n)(d_{A}/100).\\
\noalign{ \noindent{We estimate $\lceil e/2\rceil$ as $\lceil e/2 \rceil \leq \frac{C^2+1}{C^2}\frac{e}{2}$ using $e>C^2$ from Lemma \ref{alap}.
We have}}
I=\sum_{j=0}^{2} I_j &<& x_0^{2/3}Cn^{2/3}e^{2/3}+(x_1^{2/3}+x_2^{2/3})Cn^{2/3}
\left(\frac{C^2+1}{C^2} \frac{e}{2}\right)^{2/3} + \frac{(x_1+x_2)nd_A}{100} +3n+6e.\\
\noalign{ \noindent{Applying Jensen's inequality in the form $x_{1}^{2/3} +x_{2}^{2/3} \leq
2\left(\frac{x_1+x_2}{2}\right)^{2/3} =
2 \left(\frac{1-x_0}{2}\right)^{2/3}$, we obtain}}
I= \sum_{j=0}^{2} I_j &<& \left(x_0^{2/3} +\frac{(1-x_{0})^{2/3}}{2^{1/3}}
\cdot \frac{C^2+1}{C^2} + \frac{(1-x_0)}{100} \right)I
+ 6(n+e).\\
\noalign{ \noindent{By Corollary \ref{en}, we deduce}}
1&<& x_{0}^{2/3} +\frac{(1-x_{0})^{1/3}}{2^{1/3}} \cdot
\frac{(C^2+1)}{C^2}+ \frac{1-x_0}{100} +\frac{36}{C^3}.
\end{eqnarray*}
This inequality is false for $x_0\in [0, 0.1]$. (The smallest positive $x_0$
satisfying the inequality is approximately $x_0\approx 0.108$.) This proves that $x_0>0.1$,
as required
\qed

\subsection{Separation of two line sets\label{ssec:cmplx}}

Relying on the definitions of $E_1$, $E_2$, $P_0$, and $\Psi$, we formulate a lemma (Lemma~\ref{finale})
that immediately implies the Separation Lemma (Lemma~\ref{first}).

\begin{lemma}\label{finale}
Let $(P,E)$ be a critical system with $E=E_1\cup E_2$ as defined above.
Then there is a point set $P_3\subseteq P_0$ and two line sets $L_1\subseteq E_1$ and $L_2\subseteq E_2$ such that
\begin{itemize}\itemsep -2pt
\item[{\rm (a)}] $|P_3| \geq n/M^8$;
\item[{\rm (b)}] every point $p\in P_3$ is incident to at least $I/(nM^2)$ lines from each of $L_1$ and $L_2$;
\item[{\rm (c)}] there are two orthogonal directions $\hat{\ell}_1, \hat{\ell}_2\in H(1,1)$ such that an appropriate nondegenerate linear transformation of $\C_2$ maps the directions of the lines in $L_1$ and $L_2$ into the $3^{\circ}$-neighborhood of $\hat{\ell}_1$ and $\hat{\ell}_2$, respectively.
\end{itemize}
\end{lemma}

We prove Lemma~\ref{finale} at the end of this section. The main difficulty
in finding sets $L_1$ and $L_2$ with the required properties is that the
boundary $H_0$ of the two hemispheres, $H_1$ and $H_2$, is a one-dimensional manifold.
It is possible that for every point $p\in P_0$, the directions of {\em most} of the
incident lines are very close to some direction in $H_0$. This undesirable property of
a point $p\in P$ is captured in the following definition.

\begin{defini}\label{def:na}
Let $(P,E)$ be a system of points and lines in $\C^2$ such that
$E=E_1\cup E_2$ with $\widehat{E}_1\subseteq H_1$ and $\widehat{E}_2\subseteq H_2$,
and let $d_A=I/n$ denote the average number of lines in $E$ incident to a point in $P$.

For a direction $a\in H_0$, a point $p \in P$ is called an {\em $N(a)$-point},
if the directions of at least $|E_1^p|-\frac{d_A}{200\ M}$ lines in $E_1^p$ and
at least $|E_2^p|-\frac{d_A}{200\ M}$ lines in $E_2^p$ are in the open disk
of radius $10^\circ$ centered at $a\in H_0$.
\end{defini}

If a point $p\in P$ is an $N(a)$-point for some direction $a\in H_0$, then we can apply a
linear transformation $\pi_\lambda^a\in \Psi$ (cf.~Subsection~\ref{ssec:gr}) with an
appropriate parameter $\lambda\in (-1,1)$ to move all line directions out of the disk of radius
$10^\circ$. After such a transformation, however, $p$ might still be an $N(a')$-point for some
other direction $a'\in H_0$, or another point $p'\in P$ might become an $N(a')$-point
for some $a'\in H_0$. We show below (Lemma~\ref{cello}) that after applying an appropriate
linear transformation $\psi\in \Psi$, a positive fraction of the points in $P_0\subseteq P$
are no longer $N(a)$-points for any $a\in H_0$.

\begin{lemma}\label{cello}
Let $(P,E)$ be a critical system with $E=E_1\cup E_2$ as defined above.
Then there is a set $O\subseteq P_0$ of at least $n/M^6$ points and
a transformation $\psi\in \Psi$ such that after applying transformation
$\psi$, no point in $O$ is an $N(a)$-point for any $a\in H_0$.
\end{lemma}

If, after some transformation $\pi\in \Psi$,
no point in $P_0$ is an $N(a)$-point for any $a\in H_0$, then we can put $O=P_0$ in Lemma~\ref{cello}.
Otherwise, we find an appropriate set $O\subseteq P_0$ through an iterative process
over a system $(O_j,U_j\cup V_j)$ with $O_j\subseteq P_0$, $U_j\subseteq E_1$, and
$V_j\subseteq E_2$, $j\in \N_0$. Our goal is to establish Lemma~\ref{cello} for
$O=O_j$ for some $j<100$. The line sets $U_j$ and $V_j$ will be, intuitively, the
witnesses for the points in $O_j$ being $N(a)$-points for some direction $a\in H_0$.

Initially, let $O_0 = P_0$, $U_0= \{ \ell \in E_1 : \hat{\ell} \not \in H_0 \}$, and
$V_0= \{ \ell \in E_2 : \hat{\ell} \not \in H_0 \}$. That is, $U_0$ (resp., $V_0$) is obtained from $E_1$ (resp., $E_2$) by deleting all lines whose direction lies in $H_0$ (recall that this is at most one family of parallel lines). Every $p \in O_0$ is incident to at least $d_A/100-1\geq d_A/200$ lines of $U_0$ and at least $d_A/100-1\geq d_A/200$ lines of $V_0$. For $j=0$, the system $(O_j, V_j\cup U_j)$ satisfies the following four properties.

\begin{inv}\label{inv:sparse}
For $0\leq j\leq 100$, we have $O_j\subseteq P_0$, $U_j\subseteq E_1$, and $V_j\subseteq E_2$
such that
\begin{enumerate}\itemsep -2pt
\item $|O_j|  \geq n_j$, where
      $n_j = (1- \frac{3}{M})^{j} (\frac{1}{3})^j \frac{n}{10}$;
\item $| U_j \cup V_j |  \leq e_j$, where $e_j = \frac{e}{2^j}$;
\item for every $p \in O_j$, we have $|U_j^p|\geq t_j$ and  $|V_j^p|\geq t_j$,
      where $t_j= \frac{d_A}{200} (1- \frac{j}{M})$.
\end{enumerate}
\end{inv}

The following lemma describes one step of the iteration.
If we cannot choose $O$ in Lemma~\ref{cello} as a large
subset of $O_j$, then we select subsets of the points and lines
such that the number of lines decreases by a factor of 2, and the number
of lines incident to each point decreases only moderately. We show that
the process must stop after at most 100 iterations.
The number of points and lines in the system $(O_j,U_j\cup V_j)$
will monotonically decrease, but the definition of $N(a)$-points is
always understood with respect to the original system $(P,E)$,
and the original average $d_A=I/n$.

\begin{lemma}\label{iterativ}
Let $(P,E)$ be a critical system with $E=E_1\cup E_2$ as defined above.
Assume that for some $j\in \N_0$, the system $(O_j, V_j\cup U_j)$ satisfies
Invariant~\ref{inv:sparse}, and after applying any transformation $\psi\in \Psi$,
at least $(1-\frac{1}{M})n_j$ points $p\in O_j$ are each $N(a_p)$-points for some $a_p\in H_0$.
Then there are sets $O_{j+1} \subseteq O_j$,  $V_{j+1} \subset V_j$, and $U_{j+1} \subset U_j$
satisfying Invariant~\ref{inv:sparse}.
\end{lemma}

In order to prove Lemma~\ref{iterativ}, we introduce some additional notation,
and present a technical lemma (Lemma~\ref{triple}). We define two new properties
for every point $p\in P$ with respect to a set $A\subset H_0$.
\begin{quote}
For a set $A\subseteq H_0$, we say that a point $p\in P$ is an \emph{$N(A)$-point},
if $p$ is an $N(a)$-point for some $a\in H_0$ that lies in the $10^\circ$-neighborhood of $A\subset H_0$.
\end{quote}
Let $\gamma: H\rightarrow H_0$ map every direction $a\in H$, $a\not\in \{0, \infty\}$,
to the closest point on the equator $H_0$ along the main circle spanned by the directions
$a\in H$, $0\in H$ (South Pole) and $\infty\in H$ (North Pole); and let
$\gamma(0)=\gamma(\infty)=i\in H_0$.
\begin{quote}
For a set $A\subseteq H_0$, we say that a point $p\in P$ is
a {\em $\Gamma(A)$-point}, if $\gamma$ maps the directions of at least
$\frac{1}{3} |E^p|$ lines in $E^p$ to $A$.
\end{quote}

The notion of $\Gamma(A)$-points will be helpful as it is easier to track the effect 
of a linear transformation $\psi\in \Psi$ on $\Gamma(A)$-points than on $N(A)$-points.
Note that if a point $p\in \Gamma(A)$ is an $N(a)$-point for any $a\in H_0$,
then $a$ must be in the $10^{\circ}$-neighborhood of the set $A\subset H_0$.
For a set $A\subseteq H_0$ and a transformation $\psi\in \Psi$, let $N_j(A,\psi)\subseteq O_j$
denote the set of $N(A)$-points in $O_j$ after applying $\psi$. Similarly, let $\Gamma_j(A,\psi)\subseteq O_j$ denote the set of $\Gamma(A)$-points in $O_j$ after applying $\psi$.
With terminology, if a point $p\in \Gamma_j(A,\psi)$ is an $N(a)$-point for any $a\in H_0$,
then $a$ must be in the $10^{\circ}$-neighborhood of the set $A\subset H_0$, and
so $p\in N_j(A,\psi)$.

\begin{figure}[htbp]
\begin{center}
\strut\psfig{figure=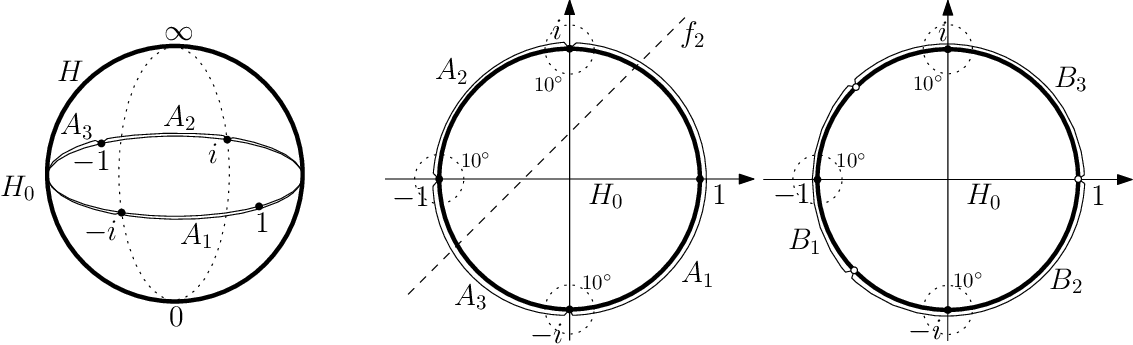,width=0.95\textwidth}
\caption{Left: a Riemann sphere $H$, the equator $H_0$ of unit slopes, and the decomposition $H_0=A_1\cup A_2\cup A_3$.
Middle: the circular arcs $A_1$, $A_2$, and $A_3$ of the equator $H_0$.
Right: the circular arcs $B_1$, $B_2$, and $B_3$ of $H_0$.\label{recurse}}
\end{center}
\end{figure}

We shall use two different decompositions of the equator $H_0$ into closed circular arcs.
First decompose the equator $H_0$ into the following three circular arcs (refer to Fig.~\ref{recurse}, left and middle): a half circle $A_1=[i,-i]$ and two quarter circles $A_2=[-i,-1]$ and $A_3=[-1,i]$.
The second decomposition of $H_0$ consists of circular arcs $H_0=B_1\cup B_2\cup B_3$
such that for every $k=1,2,3$, $A_k\cap B_k=\emptyset$ and the endpoints of $B_k$ are the midpoints
of arcs $A_{k+1 \mod 3}$ and $A_{k+2\mod 3}$.
Specifically, let $B_1=[\frac{-1+i}{\sqrt{2}},\frac{-1-i}{\sqrt{2}}]$,
$B_2=[\frac{-1-i}{\sqrt{2}},1]$, and $B_3=[1,\frac{-1+i}{\sqrt{2}}]$
(refer to Fig.~\ref{recurse}, right).
Since $H_0=A_1\cup A_2\cup A_3$ and $H_0= B_1\cup B_2\cup B_3$, every point $p\in O_j$ is a $\Gamma(A_k)$-point
for some $k\in \{1,2,3\}$, and also a $\Gamma(B_m)$-point for some $m\in \{1,2,3\}$.

\begin{lemma}\label{triple}
Let $(P,E)$ be a critical system with $E=E_1\cup E_2$ as defined above.
Assume that the system $(O_j, V_j\cup U_j)$, $j\in \N_0$, satisfies
Invariant~\ref{inv:sparse}. Then there is a transformation $\psi_0\in \Psi$ such that
$|\Gamma_j(A_k,\psi_0)|\geq \frac{n_j}{3}$ for all $k=1,2,3$.
\end{lemma}

\proof
We construct $\psi_0$ as a composition of $\pi_\lambda^1$ followed by $\pi_\kappa^i$
for some $\lambda,\kappa\in (-1,1)$. Let $\lambda\in (-1,1)$ be the minimum value
such that $|\Gamma_j(A_1,\pi_\lambda^1)|\geq n_j/3$.
This choice is possible since every point of $p\in O_j$ becomes a
$\Gamma(A_1)$-point for a sufficiently large $\lambda$, $\lambda\in (-1,1)$.
For this value of $\lambda$, we have $|\Gamma_j(A_2\cup A_3,\pi_\lambda^1)|\geq \frac{2}{3}n_j$.
In a second step, we apply $\pi_{\kappa}^i$ for some $\kappa\in (-1,1)$.
Note that for every $\kappa\in (-1,1)$, the transformation $\pi_\kappa^i$
is an automorphism on the hemisphere $\gamma^{-1}(A_1)\cup \{0,\infty\}$
and so the set of $\Gamma_j(A_1)$-points remains fixed. Let $\kappa\in (-1,1)$
be the minimum value such that $|\Gamma_j(A_2,\pi_\kappa^i\circ \pi_\lambda^1)|\geq n_j/3$
and $|\Gamma_j(A_3,\pi_\kappa^i\circ \pi_\lambda^1)|\geq n_j/3$.
This choice is possible since every point of $p\in O_j$
becomes a $\Gamma(A_3)$-point for a sufficiently large $\kappa\in (-1,1)$.
As noted above, we also have $|\Gamma_j(A_1,\pi_\lambda^1)|=
|\Gamma_j(A_1,\pi_\kappa^i\circ \pi_\lambda^1)|\geq n_j/3$.
\qed\\

We are now ready to prove Lemma~\ref{iterativ}.

\proof[Proof of Lemma~\ref{iterativ}]
Consider the system $(O_j,U_j\cup V_j)$ satisfying
Invariant~\ref{inv:sparse} such that after any transformation $\psi\in \Psi$, at
least $n_j(1-\frac{1}{M})$ points $p\in O_j$ are $N(a_p)$-points for some $a_p$.

Apply a transformation $\psi_0\in \Psi$ such that $|\Gamma_j(A_k,\psi_0)|\geq \frac{n_j}{3}$ for $k=1,2,3$
(cf.~Lemma~\ref{triple}). Recall that if a point $p\in \Gamma_j(A_k,\psi_0)$ is an $N(a_p)$-point for some
$a_p\in H_0$, then $a_p$ must be in the $10^{\circ}$-neighborhood of the arc $A_k\subset H_0$,
and so $p\in N_j(A_k,\psi_0)$. This implies that $|N_j(A_k,\psi_0)|\geq \frac{n_j}{3}-\frac{n_j}{M}=(1-\frac{3}{M})\frac{n_j}{3}$
for all $k=1,2,3$.

Since $H_0=B_1\cup B_2\cup B_3$, we have $\sum_{k=1}^3|N_j(B_k,\psi_0)|\geq n_j-\frac{n_j}{M}$ and
$\max_k|N_j(B_k,\psi_0)|\geq \frac{n_j}{3}-\frac{n_j}{3M}>(1-\frac{3}{M})\frac{n_j}{3}$.

In order to choose subsets $U_{j+1}\subseteq U_j$ and $V_{j+1}\subseteq V_j$,
we define six spherical caps of the Riemann sphere $H$, and let
$U_{j+1}\subseteq U_j$ and $V_{j+1}\subseteq V_j$ to be the set of lines whose
directions lie in one of these six spherical caps. We continue with the details.

Embed the Riemann sphere $H$ into $\R^3$ as a unit sphere $\SH^2$ centered at the origin
(i.e., every point $a\in H_0$ corresponds to a unit vector in $\R^3$). For every $a\in H_0$,
let $f(a)$ be a plane in $\R^3$ whose normal vector is $a$ and that equipartitions
the multiset of the directions $\widehat{V}_j\cup\widehat{U}_j\subset H$ (that is,
each closed halfplane bounded by $f(a)$ contains the directions of at least $|U_j\cup V_j|/2$
lines from $U_j\cup V_j$). If $a\in H$ is in general position with respect to
$\widehat{V}_j\cup\widehat{U}_j\subset H$, then $f(a)$ passes through at most
one direction in $\widehat{V}_j\cup\widehat{U}_j$. Let $a_1$, $a_2$, $a_3\in H_0$ be three
points in general position in a sufficiently small neighborhood of the midpoints of the arcs
$A_1$, $A_2$, and $A_3$, respectively (the midpoints correspond to the
directions $1\in H_0$, $(-1+i)/\sqrt{2}\in H_0$, and $(-1-i)/\sqrt{2}\in H_0$).
As a shorthand notation, let $f_1=f(a_1), f_2=f(a_2)$, and $f_3=f(a_3)$.
Refer to Fig.~\ref{recurse}. We are now ready to define the sets $O_{j+1}$, $U_{j+1}$, and $V_{j+1}$.

\begin{itemize}\itemsep -2pt
\item[-]
If there is an index $k\in \{1,2,3\}$ such that $f_k$ does not intersect the
$10^{\circ}$-neighborhood of $A_k$, then let $O_{j+1}=N_j(A_k,\psi_0)$.
Let $U_{j+1}$ (resp., $V_{j+1})$ be the set of lines from $U_j$ (resp., $V_j$) whose directions
 lie in the open spherical cap of $H$ bounded by $f_k$ that contains $B_k$.

\item[-]
If $f_k$ intersects the $10^{\circ}$-neighborhood of $A_k$ for every
$k\in \{1,2,3\}$, then consider an arc $B_m$ for $m\in \{1,2,3\}$
where $|N_j(B_m,\psi_0)|$ is maximal. Let $O_{j+1}=N_j(B_m,\psi_0)$.
Let $U_{j+1}$ (resp., $V_{j+1}$) be the set of lines from $U_j$ (resp., $V_j$)
whose directions lie in the open spherical cap of $H$ bounded by $f_m$ that contains $B_m$.
\end{itemize}

It is easy to check that $O_{j+1}$, $U_{j+1}$, and $V_{j+1}$ satisfy
Invariant~\ref{inv:sparse}. Indeed, we have $U_{j+1}\subset E_1$
and $V_{j+1}\subset E_2$ since $\psi_0\in \Psi$. We have $|O_{j+1}|\geq (1-\frac{3}{M})\frac{n_j}{3}$,
since $|N_j(A_k,\psi_0)|\geq (1-\frac{3}{M})\frac{n_j}{3}$ for $k=1,2,3$; and
$\max_k|N_j(B_k,\psi_0)|>(1-\frac{3}{M})\frac{n_j}{3}$. The number of lines is $|U_{j+1}\cup V_{j+1}|\leq e_j/2$
because of the choice of the planes $f_1$, $f_2,$ and $f_3$. Finally, each point $p\in O_{j+1}$
is an $N(a)$-point for some $a\in A_k$ or $a\in B_m$, and the set $U_{j+1}\cup V_{j+1}$ contains
all but at most $\frac{d_A}{200\ M}$ lines from $E^p$, hence also from $U_j^p\cup V_j^p$. Therefore,
we have $|U_{j+1}^p|\geq |U_j^p|-\frac{d_A}{200\ M}$ and $|V_{j+1}^p|\geq |V_j^p|-\frac{d_A}{200\ M}$, as required.
\qed

\proof[Proof of Lemma~\ref{cello}]
We count the number $I_j$ of point-line incidences for the system $(O_j,U_j)$.
On one hand, every point in $O_j$ is incident to at least $t_j$ lines in $U_j$
and so $I_j \geq |O_j|\cdot t_j$. On the other hand, the system $(O_j,U_j)$ is smaller than the
critical system $(P,E)$ and so the number of incidences is bounded above by
$\max(C|O_j|^{2/3} |U_j|^{2/3},3|O_j|,3|U_j|)$. Using $|U_j|\leq e_j$ and $|O_j|\geq n_j$, we have
\begin{eqnarray}
|O_j|\cdot t_j &\leq & \max(C|O_j|^{2/3}e_j^{2/3}, 3|O_j|,3e_j)\label{eq1}\\
t_j &\leq & \max\left(\frac{Ce_j^{2/3}}{|O_j|^{1/3}}, 3, \frac{3e_j}{|O_j|}\right).\label{eq2}\\
t_j &\leq & \max\left(\frac{Ce_j^{2/3}}{n_j^{1/3}}, 3, \frac{3e_j}{n_j}\right).\label{eq3}
\end{eqnarray}
Assuming $M=10^{10}$ and $0\leq j\leq 100$, we have $(1-j/M)\geq 1-10^{-8}$ and $0.99< (1-3/M)^j < 1$.
Lemma~\ref{alap} yields $e^{2/3}/n^{1/3}>C^2/3$ and $e^{1/3}/n^{2/3}<3/C^2$.
Consequently, we can bound the terms in~\eqref{eq3} as follows.
\begin{eqnarray}
t_j &\geq & \frac{d_A}{200}\left(1-\frac{j}{M}\right)
      > \frac{I}{203n}
      > \frac{Cn^{2/3}e^{2/3}}{203n}
      =\frac{C}{203}\cdot \frac{e^{2/3}}{n^{1/3}}
      >\frac{C^3}{609}\\
\frac{Ce_j^{2/3}}{n_j^{1/3}} & \leq & \frac{C (e/2^j)^{2/3}}{[(1-3/M)^j(1/3)^j(n/10)]^{1/3}}
                        <\frac{3C}{(4/3)^{j/3}}\cdot\frac{e^{2/3}}{n^{1/3}}\\
\frac{3e_j}{n_j}  &\leq & \frac{3e/2^j}{(1-3/M)^j(1/3)^j(n/10)}
             < \frac{31 (3/2)^j e}{n}
             < \frac{93 (3/2)^j}{C^2}\cdot \frac{e^{2/3}}{n^{1/3}}.
\end{eqnarray}
For $j=100$, we have $t_j> Ce_j^{2/3}n_j^{-1/3}$, $t_j>3$, and $t_j> 3e_jn_j^{-1}$. That is, \eqref{eq3} is false for $j=100$. Therefore there is an index $0\leq j< 100$ such that
after an appropriate transformation $\psi\in \Psi$ at least $n_j/M$ points in the
system $(O_j,U_j\cup V_j)$ are not $N(a)$-points for any $a\in H_0$. Let $O\subseteq O_j$ be
the set of these points. For every $0\leq j\leq 100$, we have
$n_j> 0.99\cdot 3^{-j} \cdot \frac{n}{10} > 3^{-100}\cdot \frac{n}{20} > n/10^{50}=n/M^5$,
hence $|O|\geq n_j/M\geq n/M^6$.
\qed\\

Cover the Riemann sphere $H$ with the minimum number of open disks (spherical caps)
of diameter $0.1^{\circ}$. Denote by $K\in \N$ the number of disks in
a minimum cover. We show that $K<M/200$ using a rough estimate.
Consider a maximal packing of the sphere $H$ with pairwise disjoint congruent
disks of radius $0.025^\circ$. The area of a disk (i.e., spherical cap) of
radius $0.025^\circ$ is more than $2(0.025 \cdot \pi/180)^2>3.807\cdot 10^{-7}$, and so the number of
disks in a packing of a unit sphere is $K<4\pi / (3.807\cdot 10^{-7}) < 3.31\cdot 10^{7}<M/200$. Increase the radius of each disk in this maximum packing from $0.025^\circ$ to
$0.05^{\circ}$ to obtain a covering of $H$ with at most $M/200$ disks. 

Partition the interior of the hemispheres $H_1$ and $H_2$ each
into at most $K$ subsets of diameter less then $0.1^{\circ}$.
Let ${\cal D}_1$ and ${\cal D}_2$, respectively, denote the families
of these subsets.

\proof[Proof of Lemma~\ref{finale}]
By Lemma~\ref{cello}, we may assume (after applying a transformation
$\psi\in \Psi$) that there is a set $O\subseteq P_0$ of $n/M^6$ points
such no point in $O$ is an $N(a)$-point for any $a\in H_0$.
Consider the partition ${\cal D}_1$ and ${\cal D}_2$ defined above.
We show that that for every point $p \in O$, we can choose two sets of 
directions, $D_1(p)\in {\cal D}_1$ and $D_2(p)\in {\cal D}_2$, such that

\begin{itemize}
\item the directions of at least $\frac{d_A}{200\cdot KM}$ lines of
      $E_1^p$ and $E_2^p$ are in $D_1(p)$ and in $D_2(p)$, respectively;
\item the distance between $D_1(p)\subset H_1$ and $D_2(p)\subset H_2$ is
      at least $6^\circ$.
\end{itemize}

Since $O\subseteq P_0$, every point $p\in O$ is incident to at
least $\frac{d_A}{200}$ lines from $E$ whose directions are in
the interior of $H_1$ (resp., $H_2$). For each $p\in O$,
choose sets $F_1(p)\in {\cal D}_1$ and $F_2(p)\in {\cal D}_2$ such
that at least $\frac{d_A}{200\cdot K}$ lines of $E_1^p$ and $E_2^p$
are in $F_1(p)$ and in $F_2(p)$, respectively.

If the distance between sets $F_1(p)$ and $F_2(p)$ is at least $6^\circ$,
then let $D_1(p)=F_1(p)$ and $D_2(p)=F_2(p)$. Otherwise let $a_p\in H_0$ be
the intersection point of the equator $H_0$ and a shortest circular arc between
$F_1(p)\subset H_1$ and $F_2(p)\subset H_2$. Assume that $F_1(p)$ is at distance
at most $3^\circ$ from $a_p$ (the case that $F_2(p)$ is at distance
at most $3^\circ$ from $a_p$ is analogous). Then $D_1(p)$ is contained in
the disk $B(a_p,3.1^{\circ})$ of radius $3.1^{\circ}$ centered at $a_p$.
Let $D_1(p)=F_1(p)$. Since $p$ is not an $N(a_p)$-point, there are at least
$\frac{d_A}{200\cdot M}$ lines in $E_2^p$ whose directions lie outside of
$B(a_p,10^{\circ})$, the disk of radius $10^\circ$ centered at $a_p$.
Out of the sets in ${\cal D}_2$ that intersect $H_2\setminus B(a_p,10^{\circ})$,
choose $D_2(p)\in {\cal D}_2$ such that it contains the directions of at least
$\frac{1}{K}\cdot \frac{d_A}{200\ M}$ lines of $E_2^p$. Since the diameter of $D_2(p)$ is
$0.1^\circ$, it lies in the exterior of $B(a_p,9.9^{\circ})$. We have
$D_1(p)\subset B(a_p,3.1^{\circ})$ and $D_2(p)\subset H_2\setminus B(a_p,9.9^{\circ})$,
and so the distance between $D_1(p)$ and $D_2(p)$ is more than $6^\circ$, as required.

For at least $|O|/ K^2$ points in $O$, we have chosen the same
subsets $D_1\in {\cal D}_1$ and $D_2\in {\cal D}_2$.
Let $P_3\subseteq O$ be the set of these points.
Since $K< M$, we have $|P_3|\geq |O|/K^2\geq (n/M^6)/K^2> n/M^8$.
Let $L_1=\{\ell\in E_1: \hat{\ell}\in D_1\}$ and
$L_2=\{\ell\in E_2: \hat{\ell}\in D_2\}$.
For every $p\in P_3$, we have $|L_1^p|\geq \frac{d_A}{200\cdot KM}>d_A/M^2$ and
$|L_2^p|\geq \frac{d_A}{200\cdot KM}>d_A/M^2$, as required.

Finally, we apply a nondegenerate linear transformation on $\C^2$ (not necessarily
from $\Psi$) that maps $D_1(p)\subset H_1$ and $D_2(p)\subset H_2$
into the $3^{\circ}$-neighborhoods of two perpendicular directions
$\hat{\ell}_1\in H$ and $\hat{\ell}_2\in H$. This can be done because
the chordal metric of $\SH^2$ is equivalent to the metric of the Riemann sphere
$H$. Let $b\in H$ be the bisector of two representative
points from $D_1(p)$ and $D_2(p)$, respectively. Apply
$\pi_\lambda^b = \varrho^b \pi_{\lambda}^1 (\varrho^b)^{-1}$,
with an appropriate $0\leq \lambda <1$, where $\varrho^b$ is an isometry
of $\C^2$ that maps the complex line of slope $b$ to a line of slope 1.
If we increase $\lambda>0$ continuously, the representative points in
$D_1(p)$ and $D_2(p)$ move along main circles of $H$ through $b\in H$.
When the representative points of the two sets become antipodal,
the diameter of the image of each set is below $3^\circ$.
\qed

\section{Covering Lemma \label{CovL}}

Our second main lemma (Covering Lemma) is an elaborate version of a
lemma of Szemer\'edi and Trotter~\cite{SzT:83b}. It states that given
a finite point set in $\R^d$, a constant fraction of the points can be
enclosed into interior-disjoint axis-aligned cubes such that the points
are approximately evenly distributed among them and each cube is adjacent
to a ``buffer zone.'' To specify what a ``buffer zone'' is, we introduce
the concepts of \emph{$\kappa$-side-cube}s and \emph{shift-graph}s.

The \emph{extents} of a hyper-rectangle $\prod_{j=1}^d[a_j,b_j]$ are the
intervals $[a_j,b_j]$, for $j=1,2,\ldots , d$. In this section, a
\emph{cube} always means an axis-aligned hypercube. A cube in $\R^d$ is
always full-dimensional unless stated otherwise. We call the direction
$\mathbf{e}_d=(0,0,\ldots, 0,1)\in \R^d$ \emph{vertical},
and interpret the ``above'' and ``below'' relationships in $\R^d$
relative to the vector ${\mathbf e}_d$.

\begin{figure}[hbp]
\begin{center}
\strut\epsfig{file=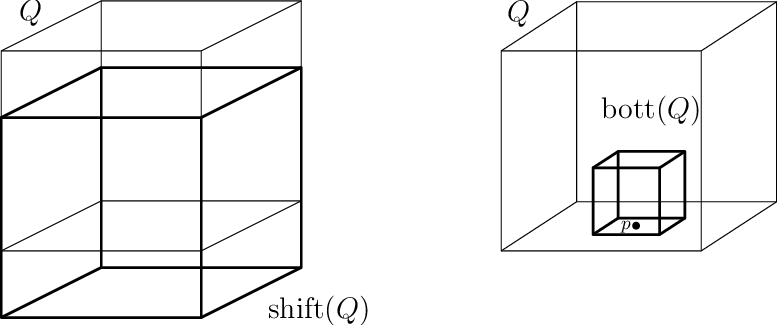,width=12cm}
\caption{Left: a cube $Q$ with ${\rm shift}(Q)$ in $\R^3$.
Right: a cube $Q$ with a $1$-side-cube ${\rm bott}(Q)$ in $\R^3$.\label{F1}}
\end{center}
\end{figure}

\begin{defini}
Let $Q$ be a cube in $\R^d$, and $\kappa \in \N$. A {\em $\kappa$-side-cube} of $Q$
is obtained by dilating $Q$ with ratio $1/(2\kappa+1)$ from a center $p$, where $p$
is the midpoint of a $(d-1)$-dimensional face of $Q$  (see Fig.~\ref{F1}).
\end{defini}

A cube has a $\kappa$-side-cube along each of its $(d-1)$-dimensional faces (sides).
So every cube in $\R^d$ has $2d$ distinct $\kappa$-side-cubes. We say that the
\emph{orientation} of a $\kappa$-side-cube $Q'$ of $Q$ is the orientation of
the vector pointing from the center of $Q$ to that of $Q'$.

\begin{defini} Let $Q$ be a cube in $\R^d$ (see Fig.~\ref{F1}).

Let $\rm{bott}(Q)$ be the $\kappa$-side-cube of $Q$ along the bottom
side of $Q$.

Let ${\rm shift}(Q)$ be the translate of $Q$ by vector
$-\frac{q}{2\kappa+1}\cdot {\mathbf e}_d$, where $q$ is
the edge length of $Q$.
\end{defini}

\begin{figure}[htbp]
\begin{center}
\strut\epsfig{file=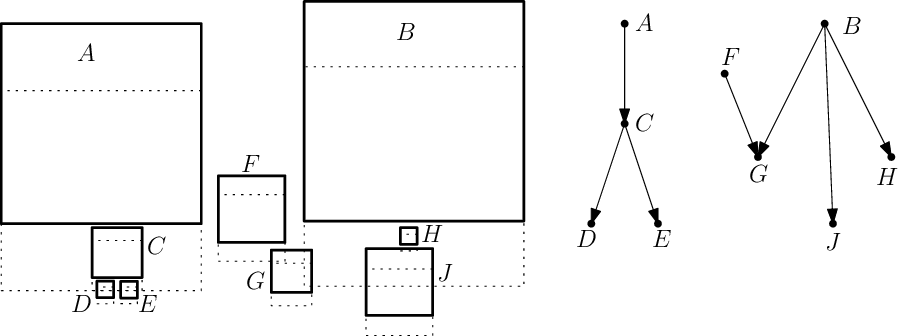,width=13cm}
\caption{\rc Left: bold interior-disjoint cubes in $\R^2$, and
dotted shifted cubes for $\kappa=1$.
Right: the corresponding shift graph.\label{TR}}
\end{center}
\end{figure}

\begin{defini}\label{def-graph}
Let ${\cal K}$ be a collection of interior-disjoint cubes in $\R^d$.
The {\em shift-graph} $T({\cal K})$ is a directed graph where
the nodes correspond to the cubes in ${\cal K}$, and
there is a directed edge $(Q_1,Q_2)$ in $T({\cal K})$
if and only if
\begin{enumerate}\topsep=0pt\itemsep=-2pt\parsep=1pt
\item  ${\rm shift}(Q_1)\setminus Q_1$ and
	${\rm shift}(Q_2)$ have a common interior point, and
\item there is a vertical segment connecting the bottom side
	of $Q_1$ and the top side of $Q_2$ that does not intersect
    the interior of any cube in ${\cal K}$. (See Fig.~\ref{TR}.)
\end{enumerate}
\end{defini}

\begin{lemma} \label{cover} {\bf (Covering Lemma)}
Let $P$ be a set of $n$ points in $\R^d$, $\kappa\in \N$, and $r\in \R$ such that
$1<r\leq n/(4(4\kappa +1)^{2d})$. Then there is a set ${\cal K}$ of
pairwise interior-disjoint cubes and a permutation of the coordinate axes
such that
\begin{enumerate}\topsep=0pt\itemsep=-2pt\parsep=1pt
\item the number of cubes is $|{\cal K}|>n/(62 d\cdot (4\kappa +1)^{2d}r)$;
\item for every cube $Q\in {\cal K}$, the interior of ${\rm bott}(Q)$
     contains at least $r$ points of $P$; and
\item the outdegree of every node in the shift graph $T({\cal K})$ is at most one.
\end{enumerate}
\end{lemma}

A permutation of the coordinate axes may be necessary, since ${\rm bott}(Q)$
and the shift graph are defined in terms of the ``vertical'' direction $\mathbf{e}_d$.

\subsection{Proof of the Covering Lemma}

A cube in $\R^d$ is called a {\em grid-cube} if all coordinates of all of its vertices are integers.
A {\em unit} grid-cube is a grid-cube of side length 1. A \emph{face-to-face tiling} of $\R^d$ with
grid-cubes is an infinite collection of pairwise interior-disjoint grid-cubes whose union is $\R^d$
such that every two grid-cubes are either disjoint or intersect in a common face.
We start with a simple proposition.

\begin{propo}\label{ppp}
Let $Q$ be a grid-cube in $\R^d$, and let $B\subseteq Q$ be a unit grid-cube.
Then $Q\setminus B$ is the union of at most $3^d-1$ (not necessarily disjoint)
grid-cubes.
\end{propo}

\proof
The hyperplanes along the $2d$ sides of $B$ decompose $Q\setminus B$ into at
most $3^d-1$ interior-disjoint axis-aligned boxes. It is enough to show
that each such box can be covered by a cube lying in $Q\setminus B$.

Consider one such box $R=\prod_{i=1}^d[a_i,b_i]\subset Q$. Note that
all coordinates of all vertices of $R$ are integers. At least one
extent of $R$ is interior-disjoint from the corresponding extent of $B$.
Assume without loss of generality that $[a_1,b_1]$ is a maximal length extent of
$R$ that is interior-disjoint from the corresponding extent of $B$.
Then a hyperplane $h_1$ orthogonal to $\mathbf{e}_1$ separates $R$ and $B$ (i.e., $B$
and $R$ lie in closed halfspaces on opposite sides of $h_1$). Enlarge each extent of
$R$ to an interval of length $b_1-a_1$ that has integer endpoints and lies in the
corresponding extent of $Q$. We obtain a cube that contains $R$, lies in $Q$, and is
separated from $Q$ by the hyperplane $h_1$.
\qed

We now state and prove a weaker version of the Covering Lemma that chooses a collection ${\cal S}$
of interior-disjoint cubes whose side-cubes jointly contain a constant fraction of the points in $P$,
but the side-cubes may have different orientations, and the condition about the shift graph is dropped.

\begin{lemma} \label{cover-}
Let $P$ be a set of $n$ points in $\R^d$, let $\kappa\in \N$, and let $r\in \R$ such that
$1<r\leq n/(4(4\kappa +1)^{2d})$. Then there is a set ${\cal S}$ of
pairwise interior-disjoint cubes such that
\begin{enumerate}\topsep=0pt\itemsep=-2pt\parsep=1pt
\item the number of cubes is $|{\cal S}|>n/(16 (4\kappa +1)^{2d}r)$;
\item for every $Q\in {\cal S}$, the interior of a $\kappa$-side-cube
      of $Q$ contains at least $r$ points of $P$.
\end{enumerate}
\end{lemma}

\proof
We are given a set of $n$ points in $\R^d$. We may assume, by applying a dilation if necessary,
that the minimum distance between any two points in $P$ is more than $\sqrt{d}$,
the diameter of a unit cube. We may also assume, by applying a translation
if necessary, that none of the coordinates of any point in $P$ is an integer.

We present a dynamic programming algorithm, Algorithm~\ref{ALG}, that
computes a collection ${\cal S}$ of cubes with the desired properties
for a point set $P\subset \mathbb{R}^d$. It proceeds in a finite number of \emph{phases}.
In phase $i\in \N$, it maintains a face-to-face tiling of $\R^d$ with a set ${\cal C}_i$ of
congruent grid-cubes. Initially, ${\cal C}_1$ is a tiling of $\R^d$
with unit grid-cubes. Our assumptions ensure that every unit grid-cube
contains at most one point from $P$, and every point in $P$ lies in the
interior of a unit grid-cube.

In each phase $i\in \N$, Algorithm~\ref{ALG} \emph{processes} every cube of ${\cal C}_i$ that
contains a point from $P$, and determines a new tiling ${\cal C}_{i+1}$.
Each cube in ${\cal C}_{i+1}$ is the union of $\mu^d$ congruent grid-cubes from ${\cal C}_i$,
where $\mu = 4\kappa +1$. Algorithm~\ref{ALG} terminates at a phase $i\in \N$, where all
points of $P$ lie in a single cube of ${\cal C}_i$. 

We allow Algorithm~\ref{ALG} to delete points from $P$, since it is enough to establish that
a side-cube of each selected cube contains \emph{at least} $r$ points of $P$.
The algorithm maintains a set $P_A\subseteq P$ of \emph{available} points.
Initially, $P_A=P$ (i.e., all points are available), but later the algorithm
may delete some points from $P_A$, or enclose them into special cubes (defined below).

When Algorithm~\ref{ALG} processes a cube $Q \in {\cal C}_i$, then it may label
$Q$ or some cubes contained in $Q$ as \emph{green}, \emph{blue}, and \emph{selected}.
Accordingly, we distinguish \emph{labeled} and \emph{unlabeled} cubes.
The green labels carry information from one phase to the next, but the blue and selected labels
are global (and irrevocable). At the end of the algorithm, the set of cubes labeled \emph{selected}
will be ${\cal S}$. We maintain the property that every two labeled cubes are either interior-disjoint
or nested. The labeled cubes are characterized as follows.

\begin{itemize}\topsep=0pt\itemsep=-2pt\parsep=1pt
\item[-]
{\bf Green cubes.} In each phase $i$, Algorithm~\ref{ALG} places some
cubes into the set ${\cal G}_i$ of \emph{green} cubes. Green cubes
are pairwise interior-disjoint. Every green cube contains at least $r$
points of $P$ but it does not contain any selected or blue cube.
At the end of phase $i+1$, some green cubes in ${\cal G}_i$
become $\kappa$-side-cubes of \emph{selected} cubes.

\item[-]
{\bf Selected cubes.} Algorithm~\ref{ALG} incrementally places cubes
into ${\cal S}$ and labels them \emph{selected}. Selected cubes are
interior-disjoint. Every $Q\in {\cal S}$ contains a green cube
as a $\kappa$-side-cube, and does not contain any smaller
blue or selected cubes.

\item[-]
{\bf Blue cubes.} Algorithm~\ref{ALG} builds a hierarchy of \emph{blue}
cubes ${\cal B}$ that enclose the selected cubes. Each blue cube contains
a unique selected cube or at least two interior-disjoint blue cubes.
\end{itemize}

In every phase $i\in \N$, Algorithm~\ref{ALG} processes all cubes in
${\cal C}_i$ that contain some points from $P$, and then determines the next tiling
${\cal C}_{i+1}$. For every tiling ${\cal C}_i$, there are $\mu^d$ possible face-to-face
tilings such that each cube in ${\cal C}_{i+1}$ contains exactly $\mu^d$ congruent cubes
from ${\cal C}_i$: Algorithm~\ref{ALG} chooses one of them to become ${\cal C}_{i+1}$.

At the end of each phase $i\in \N$, a cube $Q\in {\cal C}_i$ can be in
one of the following six states. Initially, in phase $i=1$, every cube
$Q\in {\cal C}_1$ is in state $\mathbf{A}_1$. Each state of a cube $Q\in {\cal C}_i$
is characterized by the labeled cubes contained in $Q$ that are \emph{maximal}
(for containment), and the number $|P_A\cap Q|$ of available points in $Q$
at the beginning of the step in which $Q$ is processed (note, however, that 
some of the points may be deleted from $P_A$ when $Q$ is processed).

\begin{itemize}\topsep=0pt\itemsep=-2pt\parsep=1pt
\item[$\mathbf{(A_1)}$]
	$Q$ has no label, it contains no cubes from  ${\cal G}_i\cup {\cal S}\cup {\cal B}$, and $|P_A\cap Q|<r$;

\item[$\mathbf{(A_2)}$]
  $Q\in {\cal G}$, $Q$ contains no cube from ${\cal S}\cup {\cal B}$, and $r\leq |P_A\cap Q|<\mu^d r$;

\item[$\mathbf{(A_3)}$]
 $Q$ contains one maximal blue cube and one maximal green cube $G\in {\cal G}_i$,
	and $|P_A\cap Q|<(3^d-1)r+(\mu^d-1)r < 2\mu^d r$;

\item[$\mathbf{(A_4)}$]
	$Q$ has no label, $Q$ contains one maximal blue cube, and $|P_A\cap Q|<(3^d-1)r$;

\item[$\mathbf{(A_5)}$]
	$Q\in {\cal B}$, $Q$ contains exactly one maximal selected cube,
    and $|P_A\cap Q|< \mu^d r$;

\item[$\mathbf{(A_6)}$]
	$Q\in {\cal B}$, $Q$ contains at least two maximal blue cubes,
	and $|P_A\cap Q|< 2\mu^{2d}r$;
\end{itemize}

For every cube $Q\in {\cal C}_i$ containing a point of $P$,
Algorithm~\ref{ALG} assigns $Q$ to one of the six states based on
the number $|P_A\cap Q|$ of available points in $Q$ and the states of the $\mu^d$ sub-cubes of
$Q$ from the previous subdivision ${\cal C}_{i-1}$. (By default,
all empty cubes in ${\cal C}_{i-1}$ are in state $\mathbf{A}_1$.)
We use a shorthand notation to summarize the states of all $\mu^d$
subcubes of a cube $Q\in {\cal C}_i$. The expression $Q=\sum_{k=1}^6 \omega_k \mathbf{A}_k$
means that $Q\in {\cal C}_i$ consists of $\omega_k$ sub-cubes from ${\cal C}_{i-1}$
in state $\mathbf{A}_k$, $k=1,2,\ldots ,6$ (hence $\sum_{k=1}^6 \omega_k =\mu^d$).
The assignment of a cube $Q$ to a state $\mathbf{A}_k$ is denoted by $Q\rightarrow \mathbf{A}_k$.

\begin{algorithm}\label{ALG}
Input: $P\subset \R^d$, $|P|=n$,
such that no coordinates are integers and the minimum distance between any two points is at least $\sqrt{d}$.\\
\noindent $\bullet$ Set $P_A:=P$, $i:=1$, ${\cal S}:=\emptyset$, ${\cal B}:=\emptyset$, and ${\cal G}_1:=\emptyset$.
Let ${\cal C}_1$ be the subdivision of $\R^d$ into unit grid-cubes,
each of which is in state $\mathbf{A}_1$.\\
\noindent $\bullet$ Until all points of $P_A$ lie in a single cube of ${\cal C}_i$ in state
$\mathbf{A}_1\cup \mathbf{A}_4\cup \mathbf{A}_5\cup \mathbf{A}_6$, do:

\begin{enumerate}\topsep=0pt\itemsep=-2pt\parsep=1pt
\item Set $i:=i+1$, and $G_i=\emptyset$

\item For every $Q\in {\cal C}_i$ where $Q \cap {P}\neq \emptyset$ do
\begin{enumerate}[(i)]\topsep=0pt\itemsep=-2pt\parsep=1pt
	
\item If $Q = \mu^d \mathbf{A}_1$ and $|P_A\cap Q|< r$,
       then $Q\rightarrow \mathbf{A}_1$.

\item If $Q = \mu^d \mathbf{A}_1$ and $|P_A\cap Q|\geq r$, then
	$Q\rightarrow \mathbf{A}_2$. Set ${\cal G}_i:={\cal G}_i\cup \{Q\}$
    and $P_A:=P_A\setminus Q$.

\item If $Q = (\mu^d-1) \mathbf{A}_1+ \mathbf{A}_2$, then
	$Q\rightarrow \mathbf{A}_5$. Denote by $G\subset Q$ the green subcube
   in state $\mathbf{A}_2$. $G$ is in central position within $Q$ (cf.~step~\ref{step5}).
    Let $Q^s$ be a grid-cube lying in $Q$ such that $G={\rm bott}(Q^s)$.
   Set ${\cal S}={\cal S}\cup \{Q^s\}$, ${\cal B}:={\cal B}\cup \{Q\}$, and $P_A:=P_A\setminus Q$.

\item\label{stepp++}
	If $Q = (\mu^d-1) \mathbf{A}_1+ \mathbf{A}_3$, then
	$Q\rightarrow \mathbf{A}_6$. Subcube in state $\mathbf{A}_3$ contains
    a green cube $G\in {\cal G}_{i-1}$ and a maximal blue cube $B\in {\cal B}$.
    The subcube in state $\mathbf{A}_3$ is in central position within $Q$ (cf.~step~\ref{step5}).
    Let $Q^s\subset Q\setminus B$ be a	grid-cube whose $\kappa$-side-cube is $G$.
    Set ${\cal S}:={\cal S}\cup \{Q^s\}$, ${\cal B}:={\cal B}\cup \{Q, Q^s\}$, and $P_A:=P_A\setminus Q$.

\item If $Q = (\mu^d-1)\mathbf{A}_1+(\mathbf{A}_4, \mathbf{A}_5,$  or
    $\mathbf{A}_6)$ and $|P_A\cap Q|<(3^d-1) r$, then $Q\rightarrow
	\mathbf{A}_4$.

\item\label{stepp}
	 If $Q = (\mu^d-1)\mathbf{A}_1+(\mathbf{A}_4, \mathbf{A}_5,$
	 or $ \mathbf{A}_6)$ and $|P_A\cap Q|\geq(3^d-1) r$, then
	 $Q\rightarrow \mathbf{A}_3$. Let $B$ be the maximal
     blue cube in $Q$, which either lies in the subcube
     in state $\mathbf{A}_4$ or is the subcube in state $(\mathbf{A}_5$ or $\mathbf{A}_6$).
     Let $G\subset Q\setminus B$ be one of at most $3^d-1$ cubes covering
     $Q\setminus B$ (cf.~Proposition~\ref{ppp}) such that $|P_A\cap G|\geq r$.
     Set ${\cal G}_i:={\cal G}_i\cup \{G\}$ and $P_A:=P_A\setminus G$.

\item\label{stepp+}
	If $Q = (\mu^d-2) \mathbf{A}_1+ \mathbf{A}_2 +
	(\mathbf{A}_4,\mathbf{A}_5,$ or $\mathbf{A}_6)$, then
	$Q\rightarrow \mathbf{A}_6$. Denote by $G$ and $B$ the subcubes of $Q$ in states
    $\mathbf{A}_2$ and $\mathbf{A}_4\cup \mathbf{A}_5 \cup \mathbf{A}_6$, respectively.
    The subcube in state $\mathbf{A}_2$ is in central position in $Q$ (cf.~step~\ref{step5}).
    Let $Q^s\subset Q\setminus B$ be a grid-cube whose $\kappa$-side-cube is $G$.
    Set ${\cal B}:={\cal B}\cup \{Q, Q^s\}$, ${\cal S}={\cal S}\cup	\{Q^s\}$, and $P_A:=P_A\setminus Q$.

\item
	If $Q$ contains at least two subcubes in states ($\mathbf{A}_3,
    \mathbf{A}_4,\mathbf{A}_5,$ or $\mathbf{A}_6$) and the remaining
    subcubes are in states ($\mathbf{A}_1$ or $\mathbf{A}_2$), then
	$Q\rightarrow \mathbf{A}_6$. Set ${\cal B}:={\cal B}\cup \{Q\}$
    and $P_A:=P_A\setminus Q$.
\end{enumerate}
\item Choose ${\cal C}_{i+1}$ out of the $\mu^d$ possible tilings such
that the maximal number of cubes in ${\cal C}_i$ in state $\mathbf{A}_2\cup \mathbf{A}_3$
are in central position within a cube in ${\cal C}_{i+1}$.

\item\label{step5} For every cube $Q\in {\cal C}_i$ in state
    $\mathbf{A}_2\cup \mathbf{A}_3$ that is \emph{not} in central position
    in ${\cal C}_{i+1}$, do:
    If $Q$ is in state $\mathbf{A}_2$, then
    set ${\cal G}_i:={\cal G}_i\setminus \{Q\}$ and $Q\rightarrow \mathbf{A}_1$;
    if $Q$ is in state $\mathbf{A}_3$ containing a green cube $G\in {\cal G}_i$,
    then set ${\cal G}_i:={\cal G}_i\setminus \{G\}$ and $Q \rightarrow \mathbf{A}_4$.
    As a result, every surviving cube in state $\mathbf{A}_2\cup \mathbf{A}_3$
    is in central position  within some cube in ${\cal C}_{i+1}$.

\end{enumerate}
Output: ${\cal S}$.
\end{algorithm}

At the end of Algorithm~\ref{ALG}, all points lie in a cube $Q$
in state $\mathbf{A}_1$, $\mathbf{A}_4$, $\mathbf{A}_5$ or $\mathbf{A}_6$.
Let $b$ and $s$ denote the total number of blue and selected cubes,
respectively. Let $g_i$ denote the total number of cubes that are labeled green during
phase~$i$ of Algorithm~\ref{ALG} (even if the label was removed at the end of the phase);
and let $g=\sum_{i\geq 1} g_i$.

We define a rooted tree on the blue cubes as follows.
The vertices of the tree correspond to the blue cubes,
a cube $Q_1$ is a descendant of $Q_2$ if and only if $Q_1\subset Q_2$.
Each selected cube is contained in a unique leaf cube of this tree,
and so the tree has $s$ leaves. Every intermediate node (state
$\mathbf{A}_6$) has at least two children. It follows that $b< 2s\leq 2b$.

Every green cube in ${\cal G}_i$ is in a unique cube of ${\cal C}_i$, which is in state
$\mathbf{A}_2\cup \mathbf{A}_3$. The tiling ${\cal C}_{i+1}$ was chosen such that
at least $g_i/\mu^d$ cubes $Q\in {\cal C}_i$ in $\mathbf{A}_2\cup \mathbf{A}_3$ are
in central position with respect to ${\cal C}_{i+1}$. Each of these cubes $Q\in {\cal C}_i$
contains a unique green cube $G_Q\in {\cal G}_i$. In phase~$i+1$,  a green cube
$G\in {\cal G}_i$ either becomes the $\kappa$-side-cube of a new selected cube (cases~(iii), (iv), or (vii)), or is enclosed in a new blue cube along with at least two interior-disjoint blue cubes (case~(viii)). At any rate, if a cube $Q\in {\cal C}_i$ in state $\mathbf{A}_2\cup \mathbf{A}_3$ is in
central position with respect to ${\cal C}_{i+1}$, then $Q$ is enclosed in a unique blue
cube of ${\cal C}_{i+1}$, hence $g/\mu^d\leq b\leq 2s$.

We derive an upper bound for the total number of points in terms of $s$
by accounting for the available points deleted during the algorithm.
For every blue cube (in state $\mathbf{A}_5\cup \mathbf{A}_6$), at most
$2\mu^{2d}r$ points are deleted. For every cube in state $\mathbf{A}_2\cup \mathbf{A}_3$,
at most $2\mu^dr$ points are deleted. Finally, at the last phase, the single remaining
cube is in state $\mathbf{A}_1\cup \mathbf{A}_4 \cup \mathbf{A}_5 \cup \mathbf{A}_6$,
and it contains at most $2\mu^{2d} r$ available points.
We have
\begin{eqnarray}
n&\leq & b\cdot 2\mu^{2d}r+ g\cdot 2\mu^dr+ 2\mu^{2d} r,\nonumber\\
n &<& 2s \cdot 2\mu^{2d}r + 2\mu^d s \cdot 2\mu^{d}r +2\mu^{2d} r,\nonumber\\
n-2 \mu^{2d}r &<& 8 \mu^{2d}r \cdot s,\nonumber\\
\frac{n-2 \mu^{2d} r}{8 \mu^{2d}r} &<& s.\nonumber
\end{eqnarray}
The number of cubes in ${\cal S}$ is $s > (n-2\mu^{2d} r)/(8
\mu^{2d}r) \geq n/ (16 \mu^{2d}r)$, if $r\leq n/(4 \mu^{2d})$.\qed

\proof[Proof of Lemma~\ref{cover}]
A $\kappa$-side-cube of a cube in  ${\cal S}$ can have $2d$ possible orientations.
Let ${\cal K}_0$ be the set of cubes from ${\cal S}$ with the most frequent
orientation. We can permute the coordinate axes such that every
$\kappa$-side-cube in ${\cal K}_0$ lies along the bottom sides.
The cubes of ${\cal K}_0$ satisfy properties 1 and 2 of Lemma~\ref{cover}.
The number of cubes in ${\cal K}_0$ is $|{\cal K}_0| = |{\cal S}|/2d > n/ (32d\cdot \mu^{2d}r)$.

We show that the indegree of every node in the shift graph $T({\cal K}_0)$
is at most one. Assume that $(Q_1,Q_2)$ is a directed edge in $T({\cal K}_0)$.
We wish to show that the indegree of $Q_2$ is one. Recall (Definition~\ref{def-graph})
that $(Q_1,Q_2)$ is an edge of $T({\cal K}_0)$ iff
(1) ${\rm shift}(Q_1)\setminus Q_1$ and ${\rm shift}(Q_2)$ overlap, and
(2) there is a vertical segment connecting the bottom side of $Q_1$ and
     the top side of $Q_2$ that does not intersect the interior of any cube in ${\cal K}_0$.
Since $Q_1$ and $Q_2$ are interior-disjoint, $Q_1$ must be above $Q_2$, and $Q_1$ must be
larger than $Q_2$.

It is enough to show that the vertical projection of $Q_1$ contains that of $Q_2$.
Indeed, suppose to the contrary, that both $(Q_0,Q_2)$ and $(Q_1,Q_2)$ are incoming edges
to cube $Q_2$. If the vertical projections of both $Q_0$ and $Q_1$ contain that of $Q_2$, then
one of $Q_0$ and $Q_1$ is above the other. Assume without loss of generality that $Q_0$ is above $Q_1$.
Then every vertical segment between $Q_0$ and $Q_2$ intersects $Q_1$, contradicting property (2)
of Definition~\ref{def-graph}. Therefore, the in-degree of every node in $T({\cal K}_0)$ is at most one,
as claimed.

It remains to show that the vertical projection of $Q_1$ contains that of $Q_2$.
Both $Q_1$ and $Q_2$ are in the set ${\cal S}$ of selected cubes returned
by Algorithm~\ref{ALG}. Assume without loss of generality that $Q_1$ was selected
while processing a cube $Q_1' \in {\cal C}_i$ in phase $i$. Then $Q_1'$ is the
smallest blue cube containing $Q_1$. By construction, ${\rm bott}(Q_1) \in {\cal G}_{i-1}$;
and ${\rm bott}(Q_1)$ lies in a cube of ${\cal C}_{i-1}$ in state $\mathbf{A}_2\cup \mathbf{A}_3$
that is in central position within $Q_1'$. Therefore, ${\rm shift}(Q_1)\subset Q_1'$.
Since $Q_2$ is smaller than $Q_1$ and overlaps with ${\rm shift}(Q_1)$, we have
$Q_2\subset Q_1'$. This implies that $Q_1'$ already contains some blue cube.
Hence $Q_1'$ is processed in step~\ref{stepp++} or step~\ref{stepp+} of Algorithm~\ref{ALG}.
We distinguish between two cases.

Assume first that $Q_1'$ is processed in step~\ref{stepp++}. Then $Q_1' = (\mu^d-1) \mathbf{A}_1+ \mathbf{A}_3$,
and so both ${\rm bott}(Q_1)$ and $Q_2$ lie in the cube $Q_2'$ in state $\mathbf{A}_3$, which is in central position
in $Q_1'$. Denote by $B\subset Q'_2$ the maximal blue cube in $Q_2'$.
The green cube ${\rm bott}(Q_1)\subset Q_2'$ was created in step~\ref{stepp} of Algorithm~\ref{ALG},
as one of at most $3^d-1$ cubes covering $Q\setminus B$ (cf. Proposition~\ref{ppp}).
The green cube ${\rm bott}(Q_1)$ is one of at most $3^d-1$ grid-cubes covering $Q'_2\setminus B$.
Since ${\rm bott}(Q_1)$ is adjacent to the bottom side of $Q_1$, and $Q_1$ is interior disjoint from $B$, then ${\rm bott}(Q_1)$ lies above the hyperplane containing the top side $B$.
By Proposition~\ref{ppp}, the cube ${\rm bott}(Q_1)$ is at least as large as $B$. Since
${\rm bott}(Q_1)$ and $B$ are adjacent, and ${\rm bott}(Q_1)$ is a $\kappa$-side-cube of $Q_1$ for $\kappa\geq 1$, the vertical projection of $Q_1$ contains that of $B$ and hence that of $Q_2$.

Assume now that $Q_1'$ is processed in step~\ref{stepp+}. Then $Q_1' = (\mu^d-2) \mathbf{A}_1+ \mathbf{A}_2 + (\mathbf{A}_4,\mathbf{A}_5,$ or $\mathbf{A}_6)$. The a green cube ${\rm bott}(Q_1)$ in state $\mathbf{A}_2$ is in central position in $Q_1'$, and $Q_1$ is the union of $(2\kappa+1)^d$ subcubes in ${\cal C}_{i-1}$. The cube $Q_2$ lies within a subcube $B'\in {\cal C}_{i-1}$ in state $\mathbf{A}_4\cup \mathbf{A}_5\cup \mathbf{A}_6$.
Since $(Q_1,Q_2)$ is an edge of the shift graph $T({\cal K})$, $B'$ is one of the $(2\kappa+1)^{d-1}$ subcubes in ${\cal C}_{i-1}$ directly below $Q_1$. Consequently, the vertical projection of $Q_1$ contains that of $B'$ and hence that of $Q_2$.

In both cases, the vertical projection of $Q_1$ contains the vertical projection of $Q_2$, as required.

Since the indegree of every node in the shift graph $T({\cal K}_0)$ is at most one,
it follows that $T({\cal K}_0)$ has at least $|{\cal K}_0|$ edges and so at least
half of the nodes have outdegree 0 or 1. Let ${\cal K}$ be the set of cubes in
${\cal K}_0$ whose outdegree is 0 or 1 in $T({\cal K}_0)$. We have
$|{\cal K}| \geq |{\cal K}_0|/2 >  n/ (64d\cdot \mu^{2d}r)$.
This completes the proof of Lemma~\ref{cover}.
\qed

\section{Combination of the two Main Lemmas \label{CombL}}

Recall that two 2-flats in $\R^4$ that correspond to two complex lines in $\C^2$
are either parallel or intersect in a single point. We define a \emph{crossing}
in $\R^4$ as a pair of 2-flats in $\R^4$ that intersect in exactly one point.
The Separation Lemma gives a set of points $P_3$, and two sets of 2-flats in $\R^4$,
$L_1$ and $L_2$, such that the directions of the 2-flats in $L_1$ and $L_2$ are in the
$3^\circ$-neighborhood of two orthogonal 2-dimensional subspaces
$\hat{\ell}_1\in {\rm Gr}(2,2)$ and $\hat{\ell}_2\in {\rm Gr}(2,2)$, respectively.

Lemmas~\ref{lem:localize} and \ref{lem:close} below help localizing the crossings
between the lines in $L_1$ and $L_2$. We introduce some notation for the lines incident to
two specific points. Consider two points $p,q\in \mathbb{R}^4$ and let $d={\rm dist}(p,q)$
be their Euclidean distance. Let $\ell_1^p$ and $\ell_2^p$ (reps., $\ell_1^q$ and $\ell_2^q$) be two
orthogonal 2-flats of direction $\hat{\ell}_1$ and $\hat{\ell}_2$ incident to $p$ (resp., $q$).
Let $x=\ell_1^p\cap \ell_2^q$ and $y=\ell_1^q\cap \ell_2^p$, respectively. Since $\hat{\ell}_1$ and
$\hat{\ell}_2$ are orthogonal directions, $pxqy$ is a rectangle, and in particular, $p$, $q$, $x$,
and $y$ are coplanar in $\R^4$.

\begin{lemma}\label{lem:localize}
With the above notation, the pairs of 2-flats in $L_1^p \times L_2^q$ and $L_1^q \times L_2^p$
intersect in the balls $B(x,d/10)$ and $B(y,d/10)$ of radius $\frac{d}{10}$ centered at $x$ and $y$,
respectively.
\end{lemma}
\proof
Denote by $z=h_1^p\cap h_2^q$ the intersection point of some 2-flats $h_1^p\in L_1^p$ and $h_2^q\in L_2^q$.
Since $\ell_1^p$ and $\ell_2^q$ are orthogonal, we have $\angle pxq=90^\circ$, hence $\angle pxz +\angle qxz \leq 270^\circ$.
Assume, without loss of generality, that $\angle pxz\leq 135^\circ$ (the case that $\angle qxz\leq 135^\circ$ is analogous). Recall that the distance between the directions of two 2-flats in $\R^4$
(i.e., the metric in ${\rm Gr}(2,2)$) is the sum of their principal angles. Therefore, $\angle xpz < 3^\circ$.
In the triangle $\Delta(pxz)$, we have $\angle pzx=180^\circ-\angle pxz -\angle xpz>41^\circ$. By the law of sines, ${\rm dist}(x,z)= {\rm dist}(p,x)\sin(\angle xpz)/\sin(\angle pzx)<
d \sin3^\circ/\sin 41^\circ <d/10$, as claimed.
\qed

Let $f$ be a hyperplane in $\mathbb{R}^4$, and denote its two closed halfspaces by $f^+$ and
$f^-$, respectively. A crucial step of the argument below considers the case where
two points lie on the same side of $f$, say $p,q\in f^+$, but neither $B(x,d/10)$ nor $B(y,d/10)$
is contained in $f^+$. The following lemma shows that in this case, both $p$ and $q$
must be close to $f$.

\begin{lemma}\label{lem:close}
If $p,q\in f^+$ but neither $B(x,d/10)$ nor $B(y,d/10)$ is contained in $f^+$,
then both $p$ and $q$ are at distance at most $d/5$ from the hyperplane $f$.
\end{lemma}
\proof
Since $p,q\in f^+$, the midpoint of the rectangle $pxqy$ is also in $f^+$.
This point is also the midpoint of $xy$, and so at least one of $x$ and $y$ is in $f^+$.
Assume without loss of generality that $x\in f^+$. Since $B(x,d/10)$ intersects $f^-$,
point $x$ is at distance at most $d/10$ from $f$. Hence both $x$ and $y$ are at distance
at most $d/10$ from $f$ (on either side of $f$), and the midpoint of $xy$ is also at 
distance at most $d/10$ from $f$. Since $p,q\in f^+$, and the midpoint of $pq$ is 
at distance at most $d/10$ from $f$, both $p$ and $q$ are at distance at most 
$2(d/10)=d/5$ from $f$.
\qed

We combine the Separation Lemma and the Covering Lemma in the following lemma.

\begin{lemma}{\bf (Combination Lemma)}\label{last}
Let $P$ be a set of $n$ points, and let $L_1$ and $L_2$ be two sets of 2-flats in $\R^4$
such that the directions of the 2-flats in $L_1$ and $L_2$ are in the
$3^\circ$-neighborhood of two orthogonal 2-dimensional subspaces
$\hat{\ell}_1\in {\rm Gr}(2,2)$ and $\hat{\ell}_2\in {\rm Gr}(2,2)$, respectively.
Let $r_0\in\N$, $1< r_0 \leq 10^{-8}n$.

Then there is a set ${\cal R}$ of pairwise interior-disjoint regions
in $\R^4$, and point sets $P_R \subset P\cap R$ for every $R \in {\cal R}$,
such that
\begin{enumerate}\topsep=0pt\itemsep=-2pt\parsep=1pt
\item $|{\cal R}|>n/(10^{10}r_0)$;
\item $|P_R|= r_0$ for every $R \in {\cal R}$; and
\item $\{ e_1\cap e_2: e_1\in L_1^p, e_2\in L_2^q\}\subset {\rm int}(R)$
or $\{ e_1\cap e_2: e_1\in L_1^q, e_2\in L_2^p\}\subset {\rm int}(R)$
for every $R\in {\cal R}$ and every $p,q\in P_R$.
\end{enumerate}
\end{lemma}

\proof
Invoke the Covering Lemma (Lemma~\ref{cover}) for the point set $P$
with parameters $r=27 r_0$ and $\kappa=1$ in $\R^4$. We obtain a set ${\cal K}$
of more than $n/(64\cdot 4 \cdot 5^8 \cdot 27 r_0) > 2n/(10^{10}r_0) =2n/(Mr_0)$
interior-disjoint cubes such that a 1-side-cube of each cube in ${\cal K}$
contains at least $27 r_0$ points from $P$. Assume that these special side-cubes are
the lower side-cubes ${\rm bott}(Q)$ for all $Q\in {\cal K}$
(the argument is analogous for any other orientation of the side-cubes).
For every cube $Q\in {\cal K}$, we construct a region $R\in {\cal R}$
such that $Q\cap {\rm shift}(Q) \subset R$, and choose a set $P_R$ of
$r_0$ points from $P\cap {\rm bott}(Q)$.

Consider two arbitrary points $p,q\in P\cap {\rm bott}(Q)$, for a cube $Q\in {\cal K}$.
Let $d={\rm dist}(p,q)$ and let $x$ and $y$ be defined as above.
By Lemma~\ref{lem:localize}, the pairs of lines in $L_1^p \times L_2^q$ and
$L_1^q \times L_2^p$ intersect in the balls $B(x,d/10)$ and $B(y,d/10)$ of
radius $\frac{d}{10}$ centered at $x$ and $y$, respectively. To establish
the last condition in Lemma~\ref{last}, it is enough to ensure that
at least one of the balls $B(x,d/10)$ and $B(y,d/10)$ lies in the region $R$
for all $p,q\in P_R$ (where $x$, $y$, and $d$ depend on $p$ and $q$).

The diameter of a cube in $\R^4$ is at most twice its side length, and
so $d={\rm dist}(p,q)$ is at most twice the side length of ${\rm bott}(Q)$.
It follows that both $B(x,d/10)$ and $B(y,d/10)$ are contained in ${\rm shift}(Q)$.
It is possible that for some points $p,q\in P\cap {\rm bott}(Q)$, neither $B(x,d/10)$
nor $B(y,d/10)$ is fully contained in the interior of $Q\cap {\rm shift}(Q)$. Therefore,
region $R$ should extend below $Q$.

\begin{figure}
\begin{center}
\strut\psfig{figure=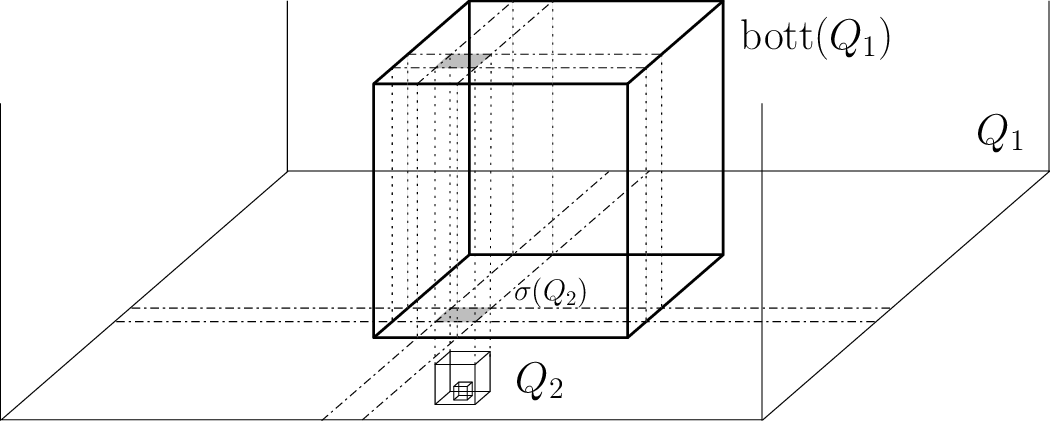,width=12cm}
\caption{Relative position of $Q_1$, ${\rm bott}(Q_1)$, and $Q_2$,
indicated in three-space instead of $\R^4$.\label{fig:cube}}
\end{center}
\end{figure}

We are now ready to define the regions $R\in {\cal R}$ and the
point sets $P_R$, $R\in {\cal R}$. We distinguish two cases.

\noindent {\bf Case 1: $Q\in {\cal K}$ has out-degree 0 in $T({\cal K})$.}
Let $R={\rm shift}(Q)$, and let $P_R$ be a set of $r_0$ arbitrary points
from $P\cap {\rm bott}(Q)$. Properties \emph{2} and \emph{3} of
Lemma~\ref{last} are satisfied for $R$, in this case ${\rm shift}(Q)$
contains both balls $B(x,d/10)$ and $B(y,d/10)$ for every $p,q\in P_R$.

\noindent {\bf Case 2: $Q_1\in {\cal K}$ has out-degree 1 in $T({\cal K})$.}
Let $Q_2\in {\cal K}$ be the cube such that $(Q_1,Q_2)$ is the unique
outgoing edge of $Q_1$. Let $\sigma:\R^4\rightarrow f_1$ be the
vertical projection of $\R^4$ to the horizontal hyperplane $f_1$;
and denote by $\sigma(Q_2)$ the vertical projection of $Q_2$.
(Fig.~\ref{fig:cube} shows the 3-dimensional analogue of the projection.)
The 2-flats spanned by the 2-dimensional faces of $\sigma(Q_2) \subset f_1$
decompose $\sigma(Q_1)$ into at most 27 axis-aligned boxes. One of them, say $Q'$,
contains the vertical projections of at least $r_0$ points in $P \cap {\rm bott}(Q_1)$.
Let $B=\{ p\in {\rm bott}(Q_1): \sigma(p)\in Q'\}$, and let $P_R$ be a set
of $r_0$ arbitrary points in $P\cap B$. Denote by ${\rm shift}(B)$ the
translate of $B$ by vector $-\frac{h}{3}\cdot {\mathbf e}_4$, where $h$ is
the side length of $Q_2$. Let region $R$ be the set of all points of
$\rm{shift}(Q_1)$ except for the points lying in or vertically below $\rm{shift}(Q_2)$.
The regions $R\in {\cal R}$ now have pairwise disjoint interiors.

It remains to show that for every $p,q\in P_R$, at least one of $B(x,d/10)$ and $B(y,d/10)$
lies in the interior of $R$. Depending on the relative position of $Q'$ and $\sigma(Q_2)$,
there are $k\in \{1,2,3,4\}$ hyperplanes that contain a 3-dimensional face of $B$ and
separate $B$ from $Q_2$. Denote these hyperplanes by $f_i$, for $i=1,\ldots, k$,
where $f_1$ is the horizontal hyperplane containing the bottom side of $Q_1$.
Denote by $f_i^+$ the halfspace bounded by $f_i$ that contains box $B$
(hence $B\subset \bigcup_{i=1}^k f_i^+$).

If $B(x,d/10)$ or $B(y,d/10)$ lies in the interior of $\bigcup_{i=1}^k f_i^+$,
then it also lies in ${\rm shift}(Q_1)\cap (\bigcup_{i=1}^k f_i^+)\subset R$.
Assume now that neither $B(x,d/10)$ nor $B(y,d/10)$ lies in the interior of $\bigcup_{i=1}^k f_i^+$.
By Lemma~\ref{lem:close}, both $p$ and $q$ are at distance at most $d/5$ from the hyperplane
$f_i$, for $i=1,\ldots , k$. The intersection of the $k$ hyperplanes, $\bigcap_{i=1}^k f_i$,
is a $(4-k)$-flat containing a face of box $B$. On one hand, the length of the orthogonal projection
of segment $pq$ to this $(4-k)$-flat is at least $\sqrt{1-k/5^2}d\geq \sqrt{1-4/5^2}d > 0.91d$.
On the other hand, the same orthogonal projection is shorter than the diameter $\sqrt{3}h$ of
the 3-dimensional cube $\sigma(Q_2)$. We have $0.91d<\sqrt{3}h$, hence $d/10<h/3$. Recall that at
least one of $x$ or $y$ is above the hyperplane $f_1$, and so the corresponding ball, $B(x,d/10)$
or $B(y,d/10)$, is strictly above the top side of $\rm{shift}(Q_2)$. It follows that $B(x,d/10)$
or $B(y,d/10)$ lies in the interior of $R$.
\qed

\subsection{Proof of Theorem~\ref{main}\label{ss-main}}

We show that the number of point-line incidences between $n$ points
and $e$ lines in $\C^2$ is at most $\max(Cn^{2/3}e^{2/3},3n,3e)$.
We proceed by contradiction. Let $(P,E)$ be a critical system of $n$
points and $e$ lines in the complex plane $\C^2$ where $n+e$ is minimal.

By the Separation Lemma, there is a set $P_3\subseteq P$ of at least
$n/M^8$ points and disjoint sets of complex lines $L_1,L_2\subset E$
such that for every point $p\in P_3$, we have $|L_1^p|\geq I/(nM^2)$ and $|L_2^p|\geq I/(nM^2)$;
and the directions of lines in $L_1$ and $L_2$ are each in a $3^{\circ}$-neighborhood
of two orthogonal directions $\hat{\ell}_1,\hat{\ell}_2\in H(1,1)$, after an appropriate
nondegenerate transformation of $\C^2$.
Identify the complex plane with the four-dimensional real Euclidean space
by $\tau: \C^2\longrightarrow \R^4$. The directions of 2-flats in $\tau(L_1)$
and $\tau(L_2)$ are each in a $3^\circ$-neighborhood of directions of
two orthogonal directions $\hat{\tau}(\hat{\ell}_2),
\hat{\tau}(\hat{\ell}_2) \in {\rm Gr}(2,2)$. For simplicity, we use the notation
$P_3=\tau(P_3)$, $L_1=\tau(L_1)$ and $L_2=\tau(L_2)$ for
a set of points and 2-flats in $\R^4$. Apply Lemma~\ref{last} for $P_3$, $L_1$, and $L_2$
with parameter $r_0=I/(nM^3)$. (The constraints $1< r_0$ and $r_0\leq 10^{-8}\cdot n/M^8$
are satisfied by Lemma~\ref{alap}.) We obtain a family ${\cal R}$ of at least
$|{\cal R}| > (n/M^8)/(Mr_0) =n/(M^9r_0)$ interior-disjoint regions in $\R^4$,
and a set $P_R\subset P_3\cap R$ of exactly $r_0$ points for each $R\in {\cal R}$.

We can now derive a lower bound for the number of crossings
$X=\{(e_1,e_2)\in L_1\times L_2\}$ between the 2-flats in
$L_1$ and $L_2$. Consider a region $R\in {\cal R}$. Each
point $p\in P_R$ is incident to at least $I/(nM^3)$ 2-flats in each of
$L_1$ and $L_2$. Since the 2-flats correspond to complex lines
and any two points determine a unique complex line, at most
$r_0$ 2-flats in $L_1^p$ (resp., $L_2^p$) may be incident to some other point
in $P_R$. Thus, there are at least $I/(nM^2) - r_0 \geq I/(2nM^2)$
2-flats in each of $L_1^p$ and $L_2^p$ that do not pass through any other
point in $P_R$. For each region $R\in {\cal R}$, we estimate
the number of crossings
$$X(P_R)=\{(e_1,e_2)\in L_1\times L_2 : \exists p,q\in P_R \mbox{ \rm such that
 }e_1\in L_1^p, e_2\in L_2^q, \mbox{ \rm and }e_1\cap e_2\in {\rm int}(R)\}.$$
By the Combination Lemma, there are at least $(I/(2nM^2))^2$ distinct crossings
for each pair $p,q\in P_R$. Every crossing is counted at most once, since
the intersection points lie in disjoint regions of ${\cal R}$. The total
number of crossings is at least
\begin{eqnarray}
|X|
&\geq & \sum_{R\in {\cal R}} |X(P_R)|
\geq |{\cal R}|\cdot {r_0\choose 2} \left(\frac{I}{2nM^2}\right)^2
> \frac{n}{M^9 r_0} \cdot \frac{r_0^2}{3} \cdot \frac{I^2}{4n^2M^4} >\nonumber\\
&>& \frac{r_0I^2}{nM^{14}}=\frac{I^3}{n^2 M^{17}} >
\frac{\max(C^3n^2e^2,27n^3,27e^3)}{n^2M^{17}}
= \frac{C^3n^2e^2}{n^2M^{17}} =\frac{C^3e^2}{M^{17}}=Me^2,\nonumber
\end{eqnarray}
by Corollary~\ref{enne} and $C^3/M^{17}=M$ (recall that $C=10^{60}$ and $M=10^{10}$).
However, $L_1, L_2 \subset E$, and so the number of crossings cannot exceed
${e\choose 2}$. The contradicting lower and upper bounds $10^{10}e^2<|X|\leq {e\choose 2}$ imply that there is no critical system $(P,E)$. We conclude that for every system of $n$ points
and $e$ lines in $\C^2$, the number of point-line incidences is bounded by
$I\leq \max(Cn^{2/3}e^{2/3},3n,3e)< Cn^{2/3}e^{2/3}+3n+3e$, as claimed.\qed

\end{document}